\documentclass[a4paper]{article}
\usepackage{amsmath,amsthm,amsfonts,amssymb}
\usepackage{multicol}
\usepackage{multirow}
\usepackage{diagbox}
\usepackage{graphicx,subfigure}

\usepackage{lineno,cases}
\allowdisplaybreaks 
\usepackage{epstopdf} 
\usepackage{pgf,tikz}
\usetikzlibrary{arrows,shapes,automata,positioning}
\usetikzlibrary{fit}
\usepackage{pst-node}
\usepackage{mathptmx}

\usepackage[pdfstartview=FitH,
CJKbookmarks=true,
bookmarksnumbered=true,
bookmarksopen=true,
colorlinks,
pdfborder=001,
linkcolor=blue,
anchorcolor=blue,
citecolor=blue
]{hyperref}

\usepackage{color,xcolor}
\renewcommand{\top}{{\mathrm{T}}}
\usepackage{bm} 
\newcommand{\lm}{\lambda}
\newcommand{\e}{\varepsilon}
\newcommand{\R}{\mathbb{R}}
\newcommand{\w}{\omega}
\newtheorem{thm}{Theorem}

\newtheorem{lem}{Lemma}
\newtheorem{ass}{Assumption}
\newtheorem{defn}{Definition}
\newtheorem{rem}{Remark}

\newcommand{\pb}{\noindent\textbf{Proof. } }
\newcommand{\pe}{\hfill\rule{4pt}{8pt}}

\def\rm{\mathrm}
\renewcommand{\top}{{\mathrm{T}}}

\begin{document}

\title{Optimal Output Consensus of High-Order Multi-Agent Systems with Embedded Technique}

\author{Yutao Tang, Zhenhua Deng, and Yiguang Hong
	\footnote{This work was supported by National Natural Science Foundation of China (61503033, 61733018, 61333001, 61573344) and National Key Research and Development Program of China (2016YFB0901902). Y. Tang is with the School of Automation, Beijing University of Posts and Telecommunications, Beijing 100876, China (e-mail: yttang@bupt.edu.cn). Z. Deng is with the School of Information Science and Engineering, Central South University, Changsha, 410004, China (e-mail: zhdeng@amss.ac.cn).  Y. Hong is with the Key Lab of Systems and Control, Institute of Systems Science, Chinese Academy of Sciences, Beijing, 100190, China, and also with the School of Mathematical Sciences, University of Chinese Academy of Sciences, Beijing 100049, China (e-mail: yghong@iss.ac.cn).}
}

\date{}

\maketitle

{\noindent\bf Abstract}: In this paper, we study an optimal output consensus problem for a multi-agent network with agents in the form of multi-input multi-output minimum-phase dynamics. Optimal output consensus can be taken as an extended version of the existing output consensus problem for higher-order agents with an optimization requirement, where the output variables of agents are driven to achieve a consensus on the optimal solution of a global cost function. To solve this problem, we first construct an optimal signal generator, and then propose an embedded control scheme by embedding the generator in the feedback loop.  We give two kinds of algorithms based on different available information along with both state feedback and output feedback, and prove that these algorithms with the embedded technique can guarantee the solvability of the problem for high-order multi-agent systems under standard assumptions.

{\noindent \bf Keywords}: Optimal output consensus, embedded control, optimal signal generator, high-order system.

\section{Introduction}

Coordination problems of multi-agent systems have
drawn much research interests due to the fast development of
large-scale systems/networks, and multiple high-order agents have
been widely discussed to deal with some practical coordination
problems, including consensus and formation in recent years
\cite{fax2004information, wang2012leader, su2012cyber, tang2015distributed}. Particularly, because of the
cost or difficulty in the measurement of all the agent states,
output consensus for high-order (minimum-phase) agents has been
widely studied \cite{rezaee2015average, seo2009consensus,xi2012output,ren2007distributed, ma2010necessary}.

In some practical applications, it might not be enough to achieve only a consensus in multi-agent systems and optimality issues of the consensus point have to be taken into consideration. Recently, distributed optimization has attracted more and more attention with their broad potential applications in multi-agent systems, smart grid and sensor networks \cite{nedic2010constrained, boyd2011distributed, bose2015equivalent,zhang2015twc}.
Among them, the case with a sum of convex objective functions has been intensively studied these years. In the problem setup, each agent is assigned with a local cost function and the control objective is to propose distributed control that guarantees a consensus on the optimal solution of the sum of all local cost functions.  Many results were obtained based on gradients or subgradients of the local cost functions combined with consensus rules. In addition to many discrete-time algorithms \cite{nedic2010constrained, yuan2011distributed, yi2014quantized, lin2016distributed}, continuous-time gradient-based optimal consensus algorithms were proposed in different situations,  considering that there are many powerful continuous-time methods. For example, in \cite{shi2013reaching}, the authors reduced the distributed optimization problem as an intersection computing problem of convex sets, where each convex set is determined by the local cost function. An alternative continuous-time algorithm was provided with discrete-time communication in \cite{kia2015distributed} to solve this problem, and moreover, a distributed optimization algorithm with disturbance rejection was studied by virtue of internal model principle in \cite{wang2015cyber}. In addition, various gradient-based algorithms have also be proposed to handle the constrained distributed optimization problem (refer to \cite{liu2015projection, yi2015distributed, liu2015second, qiu2016distributed, yang2016PI}). Note that most of the literature only consider single-integrator agents.

Furthermore, many optimal consensus or distributed optimization tasks may be implemented or depend on physical dynamics in practice, thus the distributed design involved with continuous-time agents dynamics has to be taken into account. In fact, continuous-time physical systems or hybrid cyber-physical systems are employed in practical coordination tasks, for example, in cooperative search of radio sources \cite{kim2014cooperative} and in the distributed optimal power flow problem\cite{bose2015equivalent,anese2013distributed}. Since practical systems are hardly described by single integrators, we may have to take high-order dynamics into consideration. For example, for double integrators, the authors in \cite{zhang2014distributed} proposed a distributed optimization algorithm with an integral control idea and a similar design was employed for Lagrangian agents \cite{Deng}. It was later extended to high-order integrators\cite{zhang2015distributed} assumed the existence of a kind of Lyapunov functions for high-order agents. Generally, this problem is still far from being solved.

The objective of this paper is to study the optimal output consensus problem for high-order minimum-phase multi-agent systems. To be specific, we give a framework for the optimal output consensus design over a group of multi-input multi-output (MIMO) minimum-phase agents. In light of the optimization requirement for output consensus, the gradient-based closed-loop systems are basically nonlinear since we consider general convex objective functions, although the agents are in the form of linear dynamics. The high-order features of these agents bring many new problems in its analysis and design, which make this problem much more challenging.

To overcome the difficulties brought by both high-order dynamics and nonlinearities from gradients, we
propose an embedded control scheme, consisting of an optimal signal
generator and a reference-tracking controller. This embedded technique makes
the control design carried out in a ``separative'' way, so as to simplify the
whole design by almost independently tackling the optimal consensus problem for single integrators and
tracking control of high-order agents.   In fact, the main contributions of this paper include:
\begin{itemize}
	\item[(i)] A general framework of optimal output consensus was formulated and solved for a group of MIMO minimum-phase agents with vector relative degrees. This problem can be taken as an extension by combining the two hot topics: output consensus of high-order agents \cite{ren2007distributed, ma2010necessary, xi2012output, seo2009consensus} and distributed optimization for single integrators \cite{nedic2010constrained, kia2015distributed, lin2016distributed, yi2014quantized, shi2013reaching}. Since output consensus must be achieved as the solution of a convex optimization problem of the whole multi-agent system, it is certainly more challenging than the existing output consensus problem for high-order agents. Our optimal output consensus results, in either asymptotic or exponential stability, are consistent with the existing results.
	\item[(ii)] A novel embedded control method is proposed for this optimal output consensus problem to handle the high-order dynamics in multi-agent systems.  By constructing an optimal signal generator and embedding it in the feedback loop, we can solve this problem by almost independently achieving an optimal consensus over these generators and  reference-tracking control of high-order agents.  Compared with existing optimization results for high-order integrators with some Lyapunov function assumptions in \cite{zhang2015distributed}, this embedded technique provides a quite general and constructive way to solve this problem via separating the optimization and control design.
	\item[(iii)] Two distributed gradient-based algorithms are provided by the given embedded approach in two different gradient information cases. In addition to state feedback algorithms for this problem, output feedback extensions with local observer design in each case are also proposed since the state variable of high-order agents may not be available by measurements. Under some standard assumptions, the optimal output consensus can be achieved asymptotically and even exponentially, while many existing continuous-time publications are on asymptotic stability \cite{shi2013reaching, qiu2016distributed, zhang2015distributed}.
\end{itemize}

The organization of the paper is as follows. Preliminaries about
convex analysis and graph theory are introduced and the optimal
output consensus problem is formulated in Section
\ref{sec:basis}. Then an embedded optimization control framework is
proposed for high-order minimum-phase agents with vector relative
degrees in Section \ref{sec:framework}. Main results are presented
and proved in Section \ref{sec:main} along with the given
gradient-based control in both state feedback and output feedback
cases. Following that, two examples are given to illustrate the
effectiveness of the proposed algorithm in Section \ref{sec:simu}.
Finally, concluding remarks are given in Section \ref{sec:con}.

\textsl{Notations:} Let $\R^n$ be the $n$-dimensional Euclidean
space. For a vector $x$, $||x||$ denotes its Euclidian norm. ${\bf 1}_N$ (and ${\bf 0}_N$) denotes an
$N$-dimensional all-one (and all-zero) column vector. $\mbox{col}(a_1,\,{\dots},\,a_n) =
[a_1^\top,\,{\dots},\,a_n^\top]^\top$ for column vectors $a_i\;
(i=1,\,{\dots},\,n)$. $r_N=\frac{1}{\sqrt{N}}{\bf 1}_N$, and $R_N\in \mathbb{R}^{N \times
	(N-1)}$ satisfying $R_N^\top r_N={\bf 0}_N$, $R_N^\top R_N=I_{N-1}$
and $R_N R_N^\top=I_{N}-r_N r_N^\top$.

\section{Preliminaries}\label{sec:basis}
In this section, some basic concepts are introduced for convex
analysis \cite{rockafellar1970convex, bertsekas2003convex} and graph theory \cite{godsil2001algebraic}, and
then the optimal output consensus problem is formulated.

\subsection{Convex analysis}

A function $f(\cdot)\colon \R^m \rightarrow \R $ is said to be convex if, for any $0\leq a \leq 1$,
$$ f(a\zeta_1+(1-a)\zeta_2)\leq af(\zeta_1)+(1-a)f(\zeta_2), ~ \forall \; \zeta_1,\zeta_2 \in \R^m.$$
A differentiable function $f$ is convex over $\R^m$ if
\begin{equation}\label{eq:def-convex}
f(\zeta_1)-f(\zeta_2)\geq \nabla f(\zeta_2)^T(\zeta_1 -\zeta_2),~ \forall \; \zeta_1,\zeta_2 \in \mathbb{R}^m,
\end{equation}
and $f$ is strictly convex over $\mathbb{R}^m$ if the above inequality is strict whenever $\zeta_1 \neq \zeta_2$,
and $f$ is $\omega$-strongly convex ($\omega >0$) over $\mathbb{R}^m$ if $\forall ~ \zeta_1, \zeta_2 \in \mathbb{R}^m$,
\begin{eqnarray}\label{eq:def-sc}
(\nabla f(\zeta_1)-\nabla f(\zeta_2))^T(\zeta_1 -\zeta_2)\geq \omega \|\zeta_1 -\zeta_2\|^2.
\end{eqnarray}
A vector-valued function ${\bf f}: \mathbb{R}^m \rightarrow \mathbb{R}^m$ is Lipschitz with constant $\vartheta>0$, or simply $\vartheta$-Lipschitz, if
$$
\|{\bf f}(\zeta_1)-{\bf f}(\zeta_2)\|\leq \vartheta \|\zeta_1-\zeta_2\|, ~ \forall ~ \zeta_1, \zeta_2 \in \mathbb{R}^m.
$$

\subsection{Graph theory}

A weighted undirected graph is described by $\mathcal {G}=(\mathcal
{N}, \mathcal {E}, \mathcal{A})$ with the node set
$\mathcal{N}=\{1,{\dots},N\}$ and the edge set $\mathcal {E}$
(without self-loops). $(i,\,j)\in \mathcal{E}$ denotes an edge between nodes $i$
and $j$. The weighted adjacency matrix $\mathcal{A}=[a_{ij}]\in
\mathbb{R}^{N\times N}$ is defined by $a_{ii}=0$ and
$a_{ij}=a_{ji}>0$ ($a_{ij}>0$ if and only if there is an edge
between node $i$ and node $j$). The neighbor set of node $i$ is
defined as $\mathcal{N}_i=\{j: (j,i)\in \mathcal {E} \}$ for
$i=1,\,...\,,n$.  A path in graph $\mathcal{G}$ is an alternating
sequence $i_{1}e_{1}i_{2}e_{2}{\dots}e_{k-1}i_{k}$ of nodes $i_{l}$
and edges $e_{m}=(i_{m},i_{m+1}) \in\mathcal {E}$ for
$l=1,2,{\dots},k$. If there is a path between any two vertices, then the graph is said to be connected.  The
Laplacian $L=[l_{ij}]\in \mathbb{R}^{N\times N}$ of graph
$\mathcal{G}$ is defined as $l_{ii}=\sum_{j\neq i}a_{ij}$ and
$l_{ij}=-a_{ij} (j\neq i)$, which is thus symmetric.  Denote the
eigenvalues of Laplacian matrix $L$ associated with a undirected
graph $\mathcal{G}$ as $\lambda_{1}\,\leq \,\dots\,\leq
\,\lambda_{N}$. The following lemma is well-known \cite{godsil2001algebraic}.

\begin{lem}\label{lem:graph}
	$\lambda_{1}=0$ is an eigenvalue of $L$ with ${\bf 1}_N$ as its corresponding
	eigenvector, and $\lambda_{2}>0$ if and only if the graph $\mathcal{G}$ is connected.
\end{lem}

\subsection{Problem formulation}

In this paper, we consider a group of $N$ continuous-time linear agents as follows:
\begin{align}\label{sys:agent}
\begin{cases}
\dot{\tilde{x}}_i =\tilde{A} \tilde{x}_i + \tilde B \tilde{u}_i, \\
y_i=\tilde C \tilde{x}_i, \quad i=1,\,2,\,\ldots, \,N
\end{cases}
\end{align}
where $\tilde{x}_i \in \R^{\kappa}$, $y_i \in \mathbb{R}^m$
and $\tilde u_i \in \mathbb{R}^m$ are the state variable, the output
variable and the control input, respectively.

Vector relative degree is a well-known concept and widely used in
the study of MIMO systems, which is introduced as follows.

\begin{defn}\label{def:vrd}
	System \eqref{sys:agent} is said to have a vector relative degree $(r_1,\,
	\dots, \,r_m)$ if
	\begin{itemize}
		\item[1)]  $\tilde C_i \tilde A^k \tilde B_j=0$ for all $1\leq i, \,j\leq m$ and  $0\leq k<r_i-1$
		\item[2)]  the $m\times m$ matrix  $\tilde R\triangleq [\iota_{ij}]$ is nonsingular with $\iota_{ij}=\tilde C_i \tilde A^{r_i-1} \tilde B_{j}$.
	\end{itemize}
	where $\tilde C_i$ is the $i$th row of matrix $\tilde C$ and $\tilde B_j$ is the $j$th column of $\tilde B$.
\end{defn}

To clarify our optimal output consensus problem, we focus on a class of systems as follows.

\begin{ass}\label{ass:vrd}
	System \eqref{sys:agent} is minimum-phase and has a vector relative degree $(r_1,\,...,\,r_m)$.
\end{ass}

\begin{rem}
	In fact, this assumption has been frequently used in the study of linear or nonlinear MIMO systems \cite{falb1967decoupling, isidori1995nonlinear} including integrators as special cases. Clearly, a single-input single-output linear system certainly has a relative degree.
\end{rem}

Associated with these agents, the communication topology can be described by a weighted undirected graph $\mathcal{G}=(\mathcal{V},\, \mathcal{E}, \,\mathcal{A})$, where an edge $(i,\,j)\in \mathcal{E}$ means that agent $i$ and $j$ can exchange information each other. Furthermore, each agent is endowed with a differentiable local cost function $f_i\colon \R^m \to \R$. The global cost function is defined as the sum of local costs, i.e., $f(y)=\sum_{i=1}^{N} f_i(y)$. In distributed design, agent $i$ has its own local cost function $f_i(\cdot)$, which is only known to itself and cannot be shared globally in the multi-agent network.

The optimal output consensus problem associated with those agents \eqref{sys:agent} is then formulated as follows. {\em Given a
	multi-agent system \eqref{sys:agent}, the communication graph $\mathcal{G}$ and local cost function $f_i(\cdot)$, find a
	distributed control $\tilde u_i$ for agent $i$ by only its own local data and exchanged information with its neighbors such that all the trajectories of agents are bounded and their outputs satisfy
	\begin{align}\label{def:con2}
	\lim_{t\to +\infty}||y_i(t)-y^*||=0
	\end{align}
	where $y^*$ is the optimal solution of
	\begin{align}\label{opt:main}
	\min_{ y \in \R^m} \; f(y).
	\end{align}
}

\begin{rem}\label{rem:opt-1}
	Different from the widely studied output consensus problem, optimal output consensus we formulated for high-order agents can be regarded as its extension by requiring further the consensus point as the optimal solution of a convex cost function.
	In fact, as the requirement \eqref{def:con2} coincides distributed optimization problem when agents are single integrators, it apparently implies an output consensus of the multi-agent systems \cite{ma2010necessary, xi2012output}. Many practical problems can be written in this form, including motion planing and formation control in robotics \cite{bhattacharya2011distributed, derenick2007convex}. As optimal output consensus focuses on static optimization and the steady-state performance of high-order dynamic agents, we do not give any limitations or concerns about the energy or input when we achieve the output consensus. In other words, energy/input optimization is usually involved with optimal control, but here our optimal output consensus cares about the output consensus point as the solution of a convex optimization.
\end{rem}

The following assumptions are often made to solve output consensus and/or distributed optimization problems \cite{xi2012output, kia2015distributed, wang2015cyber, nesterov2013introductory}.

\begin{ass}\label{ass:graph}
	The graph $\mathcal{G}$ is undirected and connected.
\end{ass}

\begin{ass}\label{ass:convexity-strict}
	For $i=1,\,\dots,\,N$, the function $f_i(\cdot)$ is strictly convex.
\end{ass}

Assumption \ref{ass:graph} is about the connectivity of graph $\mathcal{G}$, which guarantees that any agent's information can reach any other agents. 
Assumption \ref{ass:convexity-strict} implies the uniqueness of the optimal solution to \eqref{opt:main} \cite{rockafellar1970convex}. As usual, we assume there exists a finite optimal solution to  \eqref{opt:main} (refer to \cite{kia2015distributed, zhang2014distributed}). To get better convergence performances,  we make another assumption as follows.
\begin{ass}\label{ass:convexity-strong}
	For $i=1,\,\dots,\,N$, the function $f_i(\cdot)$ is $\w$-strongly convex and its gradient is $\vartheta$-Lipschitz for some positive constants $\w$ and $\vartheta$.
\end{ass}

Assumption \ref{ass:convexity-strong} is stronger than Assumption \ref{ass:convexity-strict}. In fact, strong convexity implies the strict convexity and the Lipschitz condition here facilitates the study of exponential convergence. Both the two assumptions have been widely used in (distributed) optimization problem \cite{nesterov2013introductory, kia2015distributed, qiu2016distributed, wang2015cyber}.

Note that the gradient-based optimization control is basically nonlinear and the general high-order system structure needs an effective design policy to optimize the output variables, which make the optimal output consensus problem much more challenging than most existing optimal consensus designs for the single-integrator and double-integrator agents.  To tackle these problems, we propose an embedded control framework in the next section.

\section{Embedded Control Scheme}\label{sec:framework}

In this section, we introduce an embedded control framework to solve the optimal output consensus problem of high-order MIMO agents in a quite unified way.

\begin{figure}
	\centering
	\newcommand\tikzscaled[2][\textwidth]{%
		\resizebox{\minof{\width}{#1}}{!}{\tikzpicture #2\endtikzpicture}}
	\pgfdeclarelayer{background}
	\pgfdeclarelayer{foreground}
	\pgfsetlayers{background,main,foreground}
	\!\!\begin{tikzpicture}[scale=0.8]
	\tikzstyle{every node}=[text centered,,font=\small,scale=0.8]
	\tikzstyle{control}=[draw, text centered, fill=blue!20,  minimum height=2em, node distance=3em]
	\tikzstyle{plant} = [control, fill=blue!40, minimum height=2.5em,text width=8em, rounded corners]
	\node (generator) [control, rounded corners, text centered, color=red!20,text=black,text width=10em] {\footnotesize Optimal signal generator};
	\node (controller) [below of=generator,text centered, node distance=6.2em,text width=14em, control, color=yellow!30,text=black] {\footnotesize \mbox{ }Reference-tracking controller};
	\node (compensator) [below of=controller,node distance=6.2em, text width=8em, minimum height= 2em, text centered, control,color=orange!20,text=black] {\footnotesize Precompensator};
	\node (plant) [below of=compensator, plant, minimum height= 2em, node distance=2.5em,  text width=8em, text=black] {\footnotesize  Agent $i$};
	\node (null) [below of=plant] {};
	\node (network)[plant, text width=5em, color=green!70, right of=null, node distance=14em, inner sep=1.3pt, text centered, text=black, minimum height=3em] {\textbf{Network}};
	\draw[thick]
	(generator.185)++(-0.8em,0)--++(-2.0em,0)--+(0,-2.2em) node[anchor=east]{$z_i$} edge[thick, ->] (controller.177)
	(controller.185)->  ++(-0.8em,0)--++(0,-3.3em)node[anchor=east]{$v_i$}edge[->]++(2em,0);
	\draw
	(compensator.185)-- ++(-0.4em , 0em) -- ++(0,-1.8em) node[anchor=east] {$\tilde u_i$} edge[->](plant.176)
	(plant.4)-- ++(0.4em,0)-- ++(0,1.8em) node[anchor=west] {$x_i$}edge[->](compensator.355);
	
	\draw[dashed,thick,->] (plant.10)++(1.8em,0)--++(2.7em,0)--++(0,7.3em)-- node[anchor=north]{$y_i$}  ++(-2.7em,0);
	\draw[thick, ->] (compensator.350)++(1.7em, 0)--++(2em,0) --++(0,1.7em) node[anchor=east]{${\chi}_i^c$}--++(0,2.1em) --(controller.0);
	
	\draw[double, ->,red]  (plant.210)++(0,-0.5em) --++(0,-2em)node[anchor=east]{$x_i$}-- ++(12em,0) --(network.195);
	\draw[double, dashed,<-,red]  (plant.330)++(0,-0.6em) --++(0,-1.00em)-- ++(6em,0) node[anchor=south]{ $x_j$}--(network.175);
	\draw[double, <-,red]  (network.125)--++(0,10em)--node[anchor=south]{$z_i$} ++(-6.7em,0);
	\draw[double, ->,red]  (network.55)--++(0,12.5em)--node[anchor=south]{$z_j$} ++(-15.5em,0)--++(0,-1.4em);
	
	\begin{pgfonlayer}{background}
	
	\path (generator.west |- generator.north)+(-0.3,0.1) node (c) {};
	\path (generator.south -| generator.east)+(0.3,-0.1) node (d) {};
	\path[fill=blue!20,rounded corners, draw=black!50,dashed]
	(c) rectangle (d);
	\path (compensator.west |- compensator.north)+(-0.6,0.1) node (c) {};
	\path (plant.south -| plant.east)+(0.6,-0.2) node (d) {};
	\path[fill=blue!20,rounded corners, draw=black!50,dashed]
	(c) rectangle (d);
	\end{pgfonlayer}
	
	\end{tikzpicture}
	\caption{Schematic diagram of the embedded control framework in optimal output consensus.}\label{fig:ctrl}
\end{figure}
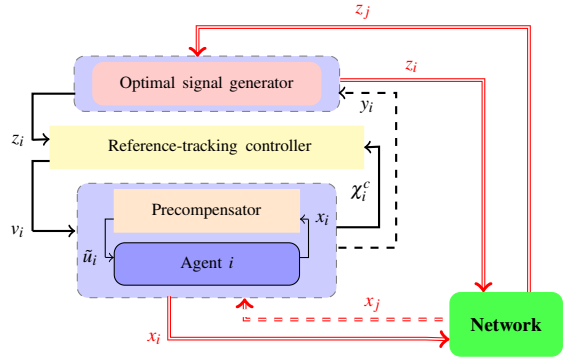

Embedded ideas or techniques are becoming ubiquitous in control applications\cite{hristu2005handbook}.   They provide flexible and reconfigurable structures to design controllers for complicated systems by integration and implementation of some standard but simple control designs. In the following, we first brief the whole idea of our embedded control scheme depicted in Fig.~\ref{fig:ctrl} and then present the first two steps of this approach.

This scheme consists of three components: {\bf precompensator} to deal with vector relative degree for \eqref{sys:agent}, {\bf optimal signal generator} to solve \eqref{opt:main}, and {\bf reference-tracking controller} to force the agent to track the generator. With a precompensator,  the original system is transformed to a normal form \eqref{sys:normal-form} with homogeneous relative degrees. This simplifies the design and makes us only focus on high-order integrators. To avoid the difficulties resulting from the high-order structure, we introduce an optimal signal generator by considering the same optimization problem for ``virtual" single integrators, in order to asymptotically reproduce the optimal solution $y^*$ by a signal $z_i$.  Then by taking $z_i$ as an output reference signal for the high-order agents and embedding this generator in the feedback loop via a tracking controller to be designed in the next section, we solve the optimal output consensus problem of system \eqref{sys:agent}.

Note that the key part for the optimization is the optimal signal generator, embedded in the feedback loop so that the physical system can achieve the optimization after its output follows the optimal signal given by this generator. In this way, we first solve the optimal (output) consensus problem for simple integrators, and then seek to synthesize the controllers for physical complex agents by embedding the optimal signal generator in a well-designed reference-tracking controller for the system \eqref{sys:normal-form} to follow \eqref{sys:n=1}.  Detailed designs of the precompensator and optimal signal generator are given in the following two subsections.

\subsection{Precompensator}

Many high-order agents have vector relative degrees with $r_i\neq r_j$ when $i\neq j$.  To make the
design simpler, we make a precondition on agents' dynamics to achieve input-output
decoupling with the same relative degree in each channel. To be specific, we transform
\eqref{sys:agent} by precompensation techniques \cite{falb1967decoupling} to a normal form through decoupling the MIMO
agent dynamics and homogenizing the relative degrees.

According to \cite{chen1995linear,isidori1995nonlinear}, we can determine a group
of basis under Assumption \ref{ass:vrd} and the associated
coordinate transformation as follows:
$$
\left\{\begin{split}
\xi^{ib}_{\iota j}= &\tilde C_\iota \tilde A^{j-1} \tilde{x}_i ~\text{with}~ j =1, \ldots, r_\iota  ~\text{and}~ \iota =1, \ldots, m, \\
\xi^{ia}_{\iota}=& \psi_\iota \tilde{x} ~\text{with}~
\psi_\iota \tilde B={\bf 0}, \iota =1, \ldots, \kappa-(r_1+\ldots+r_m).
\end{split}\right.
$$
Denote $\chi_i^b=\mbox{col}(\xi^{ib}_{1 1},\,\ldots, \,\xi^{ib}_{1
	r_1},\,\ldots, \,\xi^{ib}_{m 1},\,\ldots,\, \xi^{ib}_{m r_m})$ and
$\chi_i^a=\mbox{col}(\xi^{ia}_{1},\, \ldots,\,
\xi^{ia}_{\kappa-\sum_{j=1}^m r_j})$.  Then the system
\eqref{sys:agent} can be written in the following form:
\begin{eqnarray}\label{sys:normal-form-vrd}
\left\{\begin{split}
\dot{\chi}_i^a &= \Pi \chi_i^a +\Psi \chi_i^b, \\
\dot{\xi}^{ib}_{\iota j} &= \xi^{ib}_{\iota (j+1)}, ~ j=1,\ldots, r_\iota -1,\\
\dot{\xi}^{ib}_{\iota r_\iota} &= \Upsilon_\iota \chi_i^a +S_\iota \chi_i^b +\tilde C_\iota \tilde{A}^{r_\iota -1}\tilde B \tilde{u}_i, \\
y_i &=\mbox{col}(\xi^{ib}_{11}, \ldots, \xi^{ib}_{m1}), ~ \iota =1, \ldots, m, ~ i=1, \ldots, N,
\end{split}\right.
\end{eqnarray}
where $\Upsilon_\iota \in \mathbb{R}^{ \kappa-\sum_{j=1}^m r_j }$,
$S_\iota \in \mathbb{R}^{ \sum_{j=1}^m r_j}$, $\Psi \in
\mathbb{R}^{(\kappa-\sum_{j=1}^m r_j) \times \sum_{j=1}^m r_j}$, and
$\Pi \in \mathbb{R}^{(\kappa-\sum_{j=1}^m r_j) \times
	(\kappa-\sum_{j=1}^m r_j)}$.   Note that $\Pi$ is Hurwitz (that is,
the real parts of its eigenvalues are negative) due to the
minimum-phase of system \eqref{sys:agent}.

Take $u_i=\mbox{col}(u_{i1},\ldots, u_{im}),\; \; \Upsilon
=\mbox{col}( \Upsilon_1,\dots,       \Upsilon_m)$,
$S=\mbox{col}(S_1,\dots,S_m )$, and
\begin{equation}\label{ctrl:decoupling}
\tilde{u}_i =\tilde R^{-1} (-\Upsilon \chi_i^a -S \chi_i^b  +u_i).
\end{equation}
Therefore,
\begin{eqnarray}\label{sys:normal-form-pre}
\left\{\begin{split}
\dot{\chi}_i^a &= \Pi \chi_i^a +\Psi \chi_i^b, \\
\dot{\xi}^{ib}_{\iota j} &= \xi^{ib}_{\iota (j+1)}, ~ j=1,\ldots, r_\iota -1,\\
\dot{\xi}^{ib}_{\iota r_\iota} &= u_{i \iota}, \\
y_i &=\mbox{col}(\xi^{ib}_{11}, \ldots, \xi^{ib}_{m1}), ~ \iota =1, \ldots, m, ~ i=1, \ldots, N.
\end{split}\right.
\end{eqnarray}

Without loss of generality, we assume $\max\{r_1,..., r_m\}=n$. Then $u_i$ is taken as
\begin{equation}\label{ctrl:precompensator}
\left\{\begin{split}
u_{i \iota}^{(n-r_\iota)} &=  v_{i \iota}, ~\text{when}~ n > r_\iota,  \\
u_{i \iota} &=v_{i \iota}, ~\text{when}~ n = r_\iota,\quad \iota=1,\,\ldots, \,m
\end{split}\right.
\end{equation}
where $u_{i \iota}^{(n-r_\iota)}$ is the $(n-r_\iota)$th derivative of $u_{i \iota}$. This implies
\begin{subequations}\label{sys:normal-form}
	\begin{numcases}{}
	\dot{\chi}_i^a = \Pi \chi_i^a +\Psi \chi_i^b, \label{zds}\\
	\dot{\chi}_i^c  = (A \otimes I_m)\chi_i^c + (b \otimes I_m)v_i, \label{ssys}\\
	y_i =(c \otimes I_m) \chi_i^c, \quad i=1,\, \ldots, \,N \label{output}
	\end{numcases}
\end{subequations}
where $\chi_i^c = \mbox{col}(y_i,\ldots, y_i^{(n-1)})$, $y_i^{(j)}
=\mbox{col}( y_{i 1}^{(j)},\ldots, y_{i m}^{(j)})$,
\begin{eqnarray*}
	y_{i \iota}^{(j)} = \left\{ \begin{split}
		&  \xi^{ib}_{\iota (j+1)}, ~ j \leq r_\iota -1, \\
		&  u_{i {\iota}}, ~j=r_\iota , \\
		&  u_{i \iota}^{(j -r_\iota)} , ~ j  > r_\iota,
	\end{split}\right. ~
	\iota =1, \ldots, m, ~
	j=1, \ldots, n-1
\end{eqnarray*}
and
\begin{align*}
A=& \begin{bmatrix}
{\bf 0}_{n-1} & I_{n-1} \\
0 & {\bf 0}_{n-1}^\top \\
\end{bmatrix},  ~
b = [{\bf 0}_{n-1}^\top ~1]^\top, ~ c= [1 ~ {\bf 0}_{n-1}^\top].
\end{align*}

Consequently, we can rewrite system \eqref{sys:agent} in a normal form \eqref{sys:normal-form}, where the inputs and outputs
are decoupled and with homogeneous relative degrees $(n,\,\dots,\,n)$. In this way, we only have to consider agents in
an integrator form composed of \eqref{ssys}--\eqref{output} and get back to linear minimum-phase agents of the form \eqref{sys:agent} by a combination of \eqref{ctrl:decoupling} and well-designed $v_i$.

\subsection{Optimal Signal Generator}

Clearly, in the optimal output consensus problem, we have two main tasks: output variables of all the agents will go to the common optimal solution; and all the state variables will be bounded. Here we give an embedded technique to ``divide" the problem in order to deal with the nontrivial difficulties in the whole design. In other words, we first construct an optimal signal generator to get the optimal consensus for ``virtual" single integrators by leveraging the existing distributed optimization ideas and then use the leader-following ideas to achieve the output consensus for high-order agents, where the design for the ``virtual" single integrators can be viewed as ``embedded" technique. Here we define an optimal signal generator as follows.

\begin{defn}\label{defn:generator}
	A system $\dot{\xi}=g(\xi),\, z=h(\xi)$ is said to be an optimal
	signal generator associated with $(f_D,\,\,S,\,\,b)$ if its trajectory
	is well-defined and satisfies $z\to y^*$ as $t\to+\infty$ where
	$y^*$ solves
	\begin{align}\label{opt:defn}
	\begin{split}
	&\min_{ y \in \R^m} \quad   f_D(y)\\
	&{\rm s.t.}\quad  S^\top y=b
	\end{split}
	\end{align}
	where $S\in \R^{m\times n}$ and $b\in \R^{n}$.
\end{defn}

Obviously, an optimal signal generator is a system that can asymptotically reproduce the optimal solution of \eqref{opt:defn} (for single integrators), and thus the design of optimal signal generator is independent of available information and high-order plants. Although there may be various candidates, some of them may not lead to an augmented system for which the reference-tracking
problem can be solved. Hence, considering system composition or embedded design, we also need to concern the robustness of these generators.

Before constructing a suitable optimal signal generator for
\eqref{opt:main}, we present a useful lemma.

\begin{lem}\label{lem:key-stability}
	Consider a system of the form
	\begin{align}\label{sys:expo:lem}
	\begin{split}
	\dot{x}=-\phi(x)-Sz,\quad \dot{z}=S^\top x+T z
	\end{split}
	\end{align}
	where $x\in \R^n$,\, $z\in \R^l$, $\phi\colon \R^n \to \R^n$ is
	smooth and satisfies $\phi({\bf 0}_n)={\bf 0}_n$. Assume $T+T^\top
	\leq 0$, $(S,\, T)$ is observable, and $x^\top \phi(x) >0$ if $x\neq
	{\bf 0}_n$. Then the origin of the system is asymptotically stable.
	Furthermore, if $\phi(\cdot)$ is $\vartheta$-Lipschitz and satisfies $x^\top
	\phi(x)\geq \w x^\top x$ for some positive constants $\vartheta$ and $\w$,
	then the origin of the system is globally exponentially stable.
\end{lem}

\pb 
	Take a Lyapunov function $V=x^\top x+z^\top z$. It follows that
	\begin{align}\label{eq:lem-key}
	\dot{V}\leq -2x^\top \phi(x)+z^\top (T +T^\top)z\leq 0.
	\end{align}
	According to the LaSalle's invariance principle
	\cite{khalil2002nonlinear}, every trajectory of this system
	approaches the largest invariant set contained in ${\rm E}=\{(x,\,
	z)\mid x^\top \phi(x)=0\}$. Since no trajectory can stay identically in ${\rm E}$, other than the trivial solution $({\bf 0}_n,\, {\bf 0}_l)$ by the strict positiveness of $\phi(\cdot)$ and
	observability of the pair $(S,\,T)$.  Hence, the whole system is
	asymptotically stable at $({\bf 0}_n,\, {\bf 0}_l)$.
	
	When $\phi(\cdot)$ is $\vartheta$-Lipschitz, from its smoothness, its
	Jacobian $\frac{\partial \phi}{\partial  x}$ is globally bounded.
	Obviously, $\phi(x)=[\int_{0}^1 \frac{\partial \phi}{\partial
		x}(\theta x){\rm d}\theta ]\,x$ is well-defined. Therefore, there
	exists a time-dependent matrix $D(t)=\int_{0}^1 \frac{\partial
		\phi}{\partial x}(\theta x){\rm d}\theta $ along the trajectory of
	system \eqref{sys:expo:lem}, which is clearly uniformly bounded with
	respect to $t$. We then represent the system as follows:
	\begin{align}\label{sys:expo:lem-time-varying}
	\begin{split}
	\dot{x}=-D(t)x-Sz,\quad \dot{z}=S^\top x+T z.
	\end{split}
	\end{align}
	Note that $\bar A_o(t)=A_o-C^\top D(t) C$ with
	$$\bar A_o(t)=\begin{bmatrix}
	-D(t)& -S\\
	S^\top& T
	\end{bmatrix},\;  A_o =\begin{bmatrix}
	{\bf 0}& -S\\
	S^\top& T
	\end{bmatrix},\; C=[I_n~ {\bf 0}].$$
	Recalling the uniform boundedness of $D(t)$ and by Lemma 4.8.1 in
	\cite{ioannou1995robust} or following its proofs, we obtain that the pair $(C, \,\bar A_o(t))$
	is uniformly completely observable from the observability of $(C,\,
	A_o)$.  Hence, there exist positive constants $\delta$ and $k$ such
	that
	\begin{align*}
	W(t,\,t+\delta)=\int^{t+\delta}_t\Phi^\top(\tau,\,t)C^\top C\Phi(\tau,\, t){\rm d}\tau \geq k I_n, \quad  \forall t\geq 0
	\end{align*}
	where $\Phi(\cdot,\, \cdot)$ is the state transition matrix of \eqref{sys:expo:lem-time-varying}.
	
	The derivative of $V$ then satisfies
	\begin{align}\label{eq:lem-key-exp}
	\dot{V}\leq -x^\top \phi(x)\leq -\w x^\top x.
	\end{align}
	
	Integrating it from $t$ to $t+\delta$, we have that, for all $t\geq
	0$
	\begin{align*}
	\int^{t+\delta}_t \dot{V}(\tau){\rm d}\tau \leq - \w k (||x(t)||^2+||z(t)||^2)=-2 \w k V(t).
	\end{align*}
	Based on Theorem 8.5 in \cite{khalil2002nonlinear}, one can conclude the global
	exponential stability of this system at the origin. 
\pe 

\begin{rem}\label{rem:lem:key}
	An interesting special case is when $T=0$ and $S$ has a full-column
	rank. Apparently, $(S, \,T)$ is observable. If we further take $\phi(x)=\nabla f_D(x)$, this lemma provides the stability of the
	primal-dual gradient dynamics associated with the following
	optimization problem:
	\begin{align*}
	\begin{split}
	&\min_{x\in \R^n} ~  \quad  f_D(x)\\
	&{\rm s.t.} \quad  S^\top x={\bf 0}.
	\end{split}
	\end{align*}
\end{rem}

We then give a lemma which provides us a distributed optimal signal generator for \eqref{opt:main}.

\begin{lem}\label{lem:generator}
	Suppose Assumptions \ref{ass:graph}--\ref{ass:convexity-strict} hold. Then the optimal consensus
	problem of \eqref{opt:main} with $\dot{z}_i=v_i,\;  y_i=z_i$ is
	solved by the following control:
	\begin{align}\label{sys:n=1}
	\begin{cases}
	v_i=-\nabla f_i(z_i)-\sum\nolimits_{j=1}^N a_{ij}(\lm_i-\lm_j),\\
	\dot{\lm}_i=\sum\nolimits_{j=1}^N a_{ij}(z_i-z_j).
	\end{cases}
	\end{align}
	Furthermore, if Assumption \ref{ass:convexity-strong} is satisfied, this
	control makes $y_i(t)$ approach the optimal solution $y^*$ exponentially as
	$t\to \infty$ for $i=1,\,\ldots,\,N$.
\end{lem}

\pb  The closed-loop system can be rewritten as follows:
	\begin{align*}
	\begin{cases}
	\dot{z}=-\nabla f(z)-(L\otimes I_m)\lm,\\
	\dot{\lm}=(L\otimes I_m)z
	\end{cases}
	\end{align*}
	where $z=\mbox{col}(z_1,\,\dots,\,z_N),\,\lm=\mbox{col}(\lm_1,\,\dots,\,\lm_N)$ and $f(z)$ is determined by $f_1(z_1),\,\dots,\,f_N(z_N)$.
	
	Its equilibrium point satisfies $-\nabla f(z^*)-(L\otimes
	I_m)\lm^*={\bf 0}_N$ and $(L\otimes I_m)z^*={\bf 0}_{n\times m}$. As
	a result, there exists a $\theta$ such that
	$z_1^*=\dots=z_N^*=\theta$ since the null space of $L$ is spanned by
	${\bf 1}_N$ under Assumption \ref{ass:graph}. Then one can obtain
	$\sum_{i=1}^N \nabla f_i(\theta)=0$, which implies that $\theta$ is
	an optimal solution of \eqref{opt:main}. Thus, if we can prove the
	asymptotic stability of \eqref{sys:n=1} at its equilibrium point
	$(z^*,\,\lm^*)$, this algorithm indeed solves the distributed
	optimization problem determined by \eqref{opt:main}.
	
	For this purpose, we take $\bar z=z-z^*,\,\bar
	\lm=\lm-\lm^*$, and $\hat \lm_1=(r^\top \otimes I_m)\bar \lm,\, \hat
	\lm_2=(R^\top \otimes I_m) \bar \lm$. Since $(r^\top \otimes
	I_m)\dot{\lm}={\bf 0}_m$, it follows that $\hat \lm_1\equiv {\bf
		0}_m$ and
	$$\begin{cases}
	\dot{\bar z}=-{\bf h}-(LR\otimes I_m)\hat \lm_2,\\
	\dot{\hat \lm}_2=(R^\top L\otimes I_m)\bar z
	\end{cases}$$
	where ${\bf h}=\nabla f(z)-\nabla f(z^*)$. Denote $S=LR\otimes
	I_m$, which has a full-column rank since $R^\top L$ has a full-row rank.  Hence,
	the above error system is of the form \eqref{sys:expo:lem}.
	Note that the strict convexity of $f_i(\cdot)$ implies $\bar
	z^\top {\bf h}>0$ when $\bar z\neq{\bf 0}$. By Lemma
	\ref{lem:key-stability}, one can obtain the asymptotic stability of
	this error system and hence the solvability of this optimal consensus problem under \eqref{sys:n=1}.
	
	When $f_i(\cdot)$ satisfies Assumption \ref{ass:convexity-strong},
	it follows that $\theta=y^*$ and ${\bf h}$ is $\vartheta$-Lipschitz with $\bar
	z^\top {\bf h}\geq \w \bar z^\top \bar z$.  Then the origin of the error system is globally
	exponentially stable by Lemma \ref{lem:key-stability}, which implies that the proposed algorithm
	makes $y_i(t)$ converge to $\theta$ exponentially as $t\to \infty$ for $i=1,\,\dots,\,N$.
\pe 

\begin{rem}\label{rem:n=1}
	As we mentioned, to design an optimal signal generators is to propose a proper algorithm which can asymptotically reproduce the optimal solution of \eqref{opt:defn} (for the single integrator). In other words, the practical dynamics of high-order
	plants are not involved in its design. In the distributed design, it means to find distributed algorithms to solve the optimal consensus problem for a group of single integrators. In fact, the generator \eqref{sys:n=1} is inherently based on a primal-dual method to achieve an optimal consensus for single integrators.
\end{rem}

From Lemma \ref{lem:generator}, one can find system \eqref{sys:n=1} is an optimal signal generator for distributed optimization problem \eqref{opt:main} if $f_i(\cdot)$ is strictly convex. Furthermore, under Assumption \ref{ass:convexity-strong}, the proposed generator is of exponential stability. By Lemma 4.6 in \cite{khalil2002nonlinear}, this optimal signal generator is robust to additive perturbations on its righthand side in the sense of input-to-state stability, which will play a role in our following design. 

With the designed precompensator and optimal signal generator, we will complete the control loop by proposing proper reference-tracking controllers in the next section. In fact, this controller bridges the gap between the system \eqref{sys:normal-form} in a normal form and the optimal signal generator \eqref{sys:n=1}, and plays an important role as an interface in the design of embedded systems \cite{alur2003hierarchical}. Since $z_i$ is time-varying, set-point regulators might fail to achieve our goal. To handle this issue, we adopt some system composition techniques and a high-gain strategy in its design.

\section{Optimal Output Consensus of High-order Agents}\label{sec:main}

In this section, we complete the embedded control design and prove that the proposed control algorithms can solve the
optimal output consensus problem \eqref{opt:main} for \eqref{sys:agent} in two gradient cases.

We start with a simple case (i.e., {Case ${\bf \rm I}$}) when the gradient function $\nabla f_i(\cdot)$ is known to agent $i$, and then extend it to the second case (i.e., {Case ${\bf \rm II}$}) when only the real-time gradient $\nabla f_i(y_i)$ is available in the design. Both state and output feedback designs are proposed for two cases in the following two respective subsections.

Before the two subsections, we give the following
lemma for convergence analysis in the study of the two cases.

\begin{lem}\label{lem:cascaded}
	Consider a cascaded nonlinear system as follows:
	\begin{align}\label{sys:cascaded}
	\begin{split}
	\dot{x}_1=f_1(x_1)+g(x_1,\,x_2),\quad \dot{x}_2=f_2(x_2).
	\end{split}
	\end{align}
	Suppose for $i=1,\,2:$
	\begin{itemize}
		\item[1)] $\dot{x}_1=f_1(x_1)$ is globally exponentially stable at $x_1={\bf 0}$;
		\item[2)] $\dot{x}_2=f_2(x_2)$ is asymptotically (or globally exponentially) stable at $x_2={\bf 0}$;
		\item[3)] $||g(x_1,\,x_2)||\leq M ||x_2||$ for all $x_1,\, x_2$ with a constant $M>0$.
	\end{itemize}
	Then the system \eqref{sys:cascaded} is also asymptotically (or
	globally exponentially) stable at the origin.
\end{lem}

\pb  The asymptotic stability is a direct consequence of Lemmas 4.6 and 4.7 in \cite{khalil2002nonlinear} and thus omitted here. We only prove the global exponential stability part. From the first
	two conditions and by Theorem 4.14 in \cite{khalil2002nonlinear},
	there exist a continuously differentiable function $V_i(x)$ and
	strictly positive constants $c_{i1},\,\dots,\,c_{i4}$ such that, for
	all $x_i$,
	\begin{align*}
	&c_{i1}||x_i||^2\leq V_i(x_i)\leq c_{i2}||x_i||^2,\\
	&\dot{V}_i\leq -c_{i3}||x_i||^2, \quad ||\frac{\partial {V}_i}{\partial{x_i}}|| \leq c_{i4}||x_i||.
	\end{align*}
	
	Denoting $x\triangleq\mbox{col}(x_1,\,x_2)$ and taking a Lyapunov
	function for the cascaded system \eqref{sys:cascaded} as
	$V(x)=V_1(x_1)+cV_2(x_2)$ with $c>0$ to be determined later. Its
	derivative along \eqref{sys:cascaded} satisfies
	\begin{align*}
	\dot{V}&=\frac{\partial {V}_1}{\partial{x}_1}(f_1(x_1)+g(x_1,\,x_2))+c\frac{\partial {V}_2}{\partial{x}_2}f_2(x_2)\\
	&=-c_{13}||x_1||^2+Mc_{14}||x_1|||x_2|-cc_{23}||x_2||^2
	\end{align*}
	By Young's inequality and letting
	$c>\frac{Mc_{14}^2+2c_{13}}{2c_{13}c_{23}}$, we have
	\begin{align*}
	\dot{V}\leq -\frac{c_{13}}{2}||x_1||^2-||x_2||^2.
	\end{align*}
	Since $\min\{c_{11},\, cc_{21}\} ||x||^2\leq V(x_1,\,x_2)\leq
	\max\{c_{12},\,cc_{22}\}||x||^2$, we invoke Theorem 4.10 in
	\cite{khalil2002nonlinear} and conclude the global exponential
	stability of \eqref{sys:cascaded}.  Thus, the conclusion follows.
\pe 

\subsection{Case {I}}
In this case, the optimal signal generator can be independently implemented and asymptotically reproduce $y^*$ by Lemma \ref{lem:generator}. Thus, we only need to
propose tracking controllers for agents.

We first show our result on the high-order integrators for simplicity.

\begin{thm}\label{thm:opt-offline}
	Suppose Assumptions \ref{ass:graph}--\ref{ass:convexity-strict} hold. Then the optimal consensus
	problem \eqref{opt:main} with agents of the form ${x}_i^{(n)}=v_i,\;y_i=x_i$
	is solved by the following algorithm for any $\e>0$:
	\begin{align}\label{ctr:opt-offline}
	\begin{cases}
	v_i=-\frac{1}{\e^n}[c_0(x_i-z_i)+\sum_{k=1}^{n-1}\e^{k} c_{k} x_i^{(k)}],\\
	\dot{z}_i=-\nabla f_i(z_i)-\sum\nolimits_{j=1}^N a_{ij}(\lm_i-\lm_j),\\
	\dot{\lm}_i=\sum\nolimits_{j=1}^N a_{ij}(z_i-z_j)
	\end{cases}
	\end{align}
	where $c_0,\, \dots,\, c_{n-1}$ are constants such that
	$\sum_{k=0}^{n-1}c_{k}s^{k}+s^n$ is Hurwitz. Furthermore, if
	Assumption \ref{ass:convexity-strong} is satisfied, then $y_i(t)$ converges
	to $y^*$ exponentially as $t\to \infty$ for $i=1,\,\ldots,\,N$.
\end{thm}

\pb   Let $\hat x_i=\mbox{col}(x_i-z_i,\,\e x_i^{(1)},\,\dots,\, \e^{n-1} x_i^{(n-1)})$, it follows
	\begin{align}\label{sys:high-gain-lemma-proof-cascading}
	\dot{\hat x}_i= \frac{1}{\e} (\hat A\otimes I_m) \hat x_i+(b\otimes I_m)\dot{z}_i
	\end{align}
	where ${\small \hat A=\left[\begin{array}{c|c}
		0&I_{n-1}\\\hline
		-c_0&[\,-c_1,\,\dots,\,-c_{n-1}\,]
		\end{array}\right]}$ and $b=\mbox{col}(-1,\, {\bf 0}_{n-1})$.
	
	Then the closed-loop system is
	\begin{align}\label{sys:thm-offline-agent}
	\begin{cases}
	\dot{\hat x}_i= \frac{1}{\e} A_m \hat x_i+b_m[-\nabla f_i(z_i)-\sum\nolimits_{j=1}^N a_{ij}(\lm_i-\lm_j)],\\
	\dot{z}_i=-\nabla f_i(z_i)-\sum\nolimits_{j=1}^N a_{ij}(\lm_i-\lm_j),\\
	\dot{\lm}_i=\sum\nolimits_{j=1}^N a_{ij}(z_i-z_j)
	\end{cases}
	\end{align}
	where $A_m=\hat A\otimes I_m$ and $b_m=b\otimes I_m$.

	Due to the choice of $c_0,\, \dots,\, c_{n-1}$, the matrix $A_m$ is
	Hurwitz, which implies $||e^{(t-t_0)\frac{1}{\e} A_m}||\leq
	k_1e^{-\frac{k_2}{\e}(t-t_0)}$ for some positive constants $k_1$ and
	$k_2$. Solving \eqref{sys:high-gain-lemma-proof-cascading} and using this
	bound, we obtain
	\begin{align*}
	||\hat x(t)||
	&\leq k_1e^{-\frac{k_2}{\e}(t-t_0)}||\hat x(t_0)||+\frac{k_1\e||b_m||}{k_2}\sup_{t_0\leq \tau \leq t}||\dot{z}(\tau)||.
	\end{align*}
	By Lemma \ref{lem:generator}, $z_i$ converges to $y^*$ and
	$\dot{z}_i$ converges to $\bf 0$ when $t\to \infty$ for $i=1,\,\ldots,\,N$. For any given $\delta>0$, there is a time $T_1>0$
	such that $||\dot{z}(t)||\leq \frac{\delta k_2}{2 k_1\e||b_m||}$ for
	all $t>T_1$. Choose $t_0>T_1$ and there exists a $T_2>0$ such that
	$k_1e^{-\frac{k_2}{\e}(t-t_0)}||x(t_0)||\leq \frac{\delta}{2}$ for
	$t>T_2$. Thus $||\hat x(t)||\leq \delta$ for all
	$t>\max\{T_1,\,T_2\}$, which leads to the convergence of $x_i-z_i$
	for $i=1,\,\ldots,\,N$. Recalling $\lim_{t\to+\infty}z_i=y^*$ gives the convergence of $x_i$ with respect to $y^*$ for $i=1,\,\ldots,\,N$.
	
	To prove the exponential convergence, we rewrite the whole closed-loop system as follows:
	\begin{align}\label{sys:thm:opt-offline-proof}
	\begin{cases}
	\dot{\hat x}_i= \frac{1}{\e} A_m \hat x_i+b_m\dot{\bar z}_i,\\
	\dot{\bar z}_i=-{\bf h}_i-\sum\nolimits_{j=1}^N a_{ij}(\bar\lm_i-\bar \lm_j),\\
	\dot{\bar \lm}_i=\sum\nolimits_{j=1}^N a_{ij}(\bar z_i-\bar z_j)
	\end{cases}
	\end{align}
	where $\hat x_i=\mbox{col}(x_i-z_i,\,\e x_i^{(1)},\,\dots,\,\e^{n-1}
	x_i^{(n-1)})$, $\bar z_i=z_i-y^*,\, \bar \lm_i=\lm_i-\lm_i^*$, and
	${\bf h}_i=\nabla f_i(z_i)-\nabla f_i(y^*)$ for $i=1,\,\ldots,\,N$.
	
	If Assumption \ref{ass:convexity-strong} holds, $(z_i,\,\lm_i)$
	exponentially converges to its equilibrium point $(y^*,\,\lm_i^*)$
	as $t\to \infty$ for $i=1,\,\ldots,\,N$ by Lemma \ref{lem:generator}. From
	the Lipschitzness of $\nabla f_i(\cdot)$, the cross-term $\dot{\bar
		z}_i$ is also Lipschitz with respect to $(\bar z_i,\,\bar \lm_i)$.
	By Lemma \ref{lem:cascaded}, the exponential stability of
	\eqref{sys:thm:opt-offline-proof} is achieved, which implies the
	conclusion.
\pe 

\begin{rem}\label{rem:strict}
	Clearly, as the gradients of local cost functions are analytically known, the optimal signal generator can be implemented independently. Then the whole system is in a cascading form and the  strict convexity will suffice the solvability of this optimal output consensus problem. When Assumption \ref{ass:convexity-strong} is satisfied, the convergence can be exponentially fast by Lemma \ref{lem:cascaded}.
\end{rem}

Based on Theorem \ref{thm:opt-offline} and the normal form \eqref{sys:normal-form}, we go back to investigate linear
minimum-phase agents of the form \eqref{sys:agent} and get the following controller:
\begin{equation}\label{ctr:linear-offline}
\begin{cases}
\tilde{u}_i=\tilde R^{-1} (-\Upsilon \chi_i^a -S \chi_i^b +u_i),\\
v_i=-\frac{1}{\e^n}[c_0(y_i-z_i)+\sum_{k=1}^{n-1}\e^{k} c_{k} y_i^{(k)}],\\
\dot{z}_i=-\nabla f_i(z_i)-\sum\nolimits_{j=1}^N a_{ij}(\lm_i-\lm_j),\\
\dot{\lm}_i=\sum\nolimits_{j=1}^N a_{ij}(z_i-z_j)
\end{cases}
\end{equation}
where $c_1, \dots, c_{n-1}$ are selected as given in Theorem
\ref{thm:opt-offline}, and the relationship between $u_i$ and $v_i$
is described by \eqref{ctrl:precompensator}.

\begin{thm}\label{thm:opt-offline-linear}
	Suppose Assumptions \ref{ass:vrd}-\ref{ass:convexity-strict} hold. Then the
	optimal output consensus problem of system \eqref{sys:agent} is
	solved by  \eqref{ctr:linear-offline} for any $\e>0$.  Furthermore,
	if Assumption \ref{ass:convexity-strong} is satisfied, $y_i(t)$ converges
	to $y^*$ exponentially as $t\to \infty$ for $i=1,...,N$.
\end{thm}

\pb  Based on the above analysis, the output variables of system \eqref{sys:normal-form} under control $v_i$ and system
	\eqref{sys:agent} under control \eqref{ctr:linear-offline} are the same. Thus, we only have to investigate the convergence of
	\eqref{sys:normal-form} under the above control.  Clearly, the subsystem \eqref{ssys}-\eqref{output} is equivalent to
	the $n$th-order integrator.  Let ${\tilde \chi}_i^a={\chi}_i^a-{{\chi}_i^a}^*, \,{\tilde \chi}_i^b={\chi}_i^b-{{\chi}_i^b}^*$ (where ${{\chi}_i^a}^*, \,{{\chi}_i^b}^*$ are the equilibrium point of 	\eqref{zds}-\eqref{ssys} determined by $y^*$) and define ${\hat \chi}^c_i$ as $\hat x_i$ in Theorem \ref{thm:opt-offline}. The whole systems is rewritten as:
	$$
	\begin{cases}
	\dot{\tilde \chi}_i^a = \Pi \tilde\chi_i^a +\Psi \tilde\chi_i^b,\\
	\dot{\hat \chi}_i^c= \frac{1}{\e} A_m {\hat \chi}_i^c+b_m\dot{\bar z}_i,\\
	\dot{\bar z}_i=-{\bf h}_i-\sum\nolimits_{j=1}^N a_{ij}(\bar \lm_i-\bar \lm_j),\\
	\dot{\bar \lm}_i=\sum\nolimits_{j=1}^N a_{ij}(\bar z_i-\bar z_j).\\
	\end{cases}
	$$
	
	Apparently,  it is in a cascading form with the first subsystem as a driven one. Recalling the minimum-phase property of \eqref{sys:agent}, the
	matrix $\Pi$ is Hurwitz. According to Theorem \ref{thm:opt-offline}, the output $y_i$ of the last three subsystems asymptotically converges to the
	optimal solution $y^*$, then the asymptotic convergence of the whole cascaded system can be obtained by Lemma \ref{lem:cascaded}. When Assumption \ref{ass:convexity-strong} holds, the global exponential stability 	can be obtained in a similar way. The proof is thus complete.
\pe

Next, let us consider the case when only the output variables of each agent can be obtained because it may be difficult to get or
measure all the state variables in some situations.  Since the optimal signal generator is independently implemented, we only have
to focus on the tracking part. To solve the problem, we consider an output feedback version of the proposed
high-gain embedded control by proposing an observer-based output feedback design.

Since system \eqref{sys:agent} is minimum phase and therefore detectable, we can design the following local
observer for agent $i$ ($i =1, \ldots, N$) of the form
\eqref{sys:agent}.
\begin{align}
\label{eq.mpobserver}
\begin{cases}
\dot{\hat{y}}_i=\hat{y}_i^{(1)} - l_1 (\hat{y}_i-y_i), \\
\dot{\hat{y}}_i^{(1)}=\hat{y}_i^{(2)}- l_2 (\hat{y}_i-y_i),\\
\cdots\\
\dot{\hat{y}}_i^{(n-1)}=v_i- l_n (\hat{y}_i-y_i),\\
\dot{\hat{\chi}}_i^a = \Pi \hat{\chi}_i^a +\Psi \hat\chi_i^b\end{cases}
\end{align}
where $\hat{\chi}_i^a$, $\hat{\chi}_i^b$, $\hat{y}_i$ and
$\hat{y}_i^{(\iota)}$ are the estimations of ${\chi}_i^a$,
$\chi_i^b$, $y_i$ and ${y}_i^{(\iota)}$, respectively, with $\iota
\in \{1, \ldots, n-1\}$, and $l_1, \ldots, l_n$ are constants such
that $p_l(s)=s^n + l_1s^{n-1}+ \ldots + l_{n-1}s +l_n$ is Hurwitz.
Then we substitute these variables by their estimations and propose
the following distributed control for system \eqref{sys:agent}.
\begin{equation}\label{ctr:opt-offline-output-linear}
\begin{cases}
\tilde{u}_i=\tilde R^{-1} (-\Upsilon \hat \chi_i^a -S \hat\chi_i^b +u_i),\\
v_i=-\frac{1}{\e^n}[c_0(y_i-z_i)+\sum_{k=1}^{n-1}\e^{k} c_{k} \hat y_i^{(k)}],\\
\dot{z}_i=-\nabla f_i(z_i)-\sum\nolimits_{j=1}^N a_{ij}(\lm_i-\lm_j),\\
\dot{\lm}_i=\sum\nolimits_{j=1}^N a_{ij}(z_i-z_j)\\
\end{cases}
\end{equation}
where $c_0, \, \ldots, \, c_{n-1}$ and  $\e>0$  are defined as before, the
relationship between $u_i$ and $v_i$ is described by \eqref{ctrl:precompensator}.

We then give the following results for the output feedback design.
\begin{thm}\label{thm:offline-ouput-linear}
	Suppose Assumptions \ref{ass:vrd}--\ref{ass:convexity-strict} hold. Then the
	optimal output consensus problem for the agents of the form
	\eqref{sys:agent} is solved by \eqref{ctr:opt-offline-output-linear}
	for any $\e>0$. Furthermore, if Assumption
	\ref{ass:convexity-strong}  is satisfied, $y_i(t)$ converges to $y^*$
	exponentially as $t\to \infty$ for $i=1,\,\ldots,\,N$.
\end{thm}
\pb    Taking $\bar e_i=\mbox{col}(y_i-\hat y_i,\dots,
	y_i^{(n)}-\hat{y}_i^{(n)})$ and $\bar  \chi_i^a=\hat{\chi}_i^a-{\chi}_i^a,\,
	\bar{\chi}_i^b=\hat{\chi}_i^b-{\chi}_i^b$ gives
	$$\begin{cases}
	\dot{\bar e}_i=(A_l\otimes I_m)\bar e_i,\\
	\dot{\bar{\chi}}_i^a = \Pi \bar{\chi}_i^a +\Psi \bar\chi_i^b\end{cases}$$
	where $$ A_l=\begin{bmatrix}-{\bf l}_{n-1}&I_{n-1}\\-l_n&{\bf 0}^\top_{n-1}\end{bmatrix},\quad {\bf l}_{n-1}=\mbox{col}(l_1,\dots,l_{n-1}).$$
	
	Recalling the Hurwitzness of $\Pi$, we apply Lemma \ref{lem:cascaded} and obtain $\hat{\chi}_i^a$, $\hat{\chi}_i^b$, $\hat{y}_i$, $\hat{y}_i^{(\iota)}$ exponentially converge to ${\chi}_i^a$, $\chi_i^b$, ${y}_i$ and
	${y}_i^{(\iota)}$ with $\iota \in \{1,\, \ldots,\, n-1\}$ for $i \in  \{1,2, \ldots, N\}$, respectively.
	
	Since substituting ${\chi}_i^a$, $\chi_i^b$, ${y}_i$ and ${y}_i^{(\iota)}$ by their estimations will not change the equilibrium point, we only have to prove the stability of the new system with respect to its equilibrium point under observer-based control laws. For this purpose, we repeat the whole error systems as follows:
	$$ 	
	\begin{cases}
	\dot{\tilde \chi}_i^a = \Pi \tilde\chi_i^a +\Psi \tilde\chi_i^b,\\
	\dot{\hat \chi}_i^c= \frac{1}{\e} A_m {\hat \chi}_i^c+b_m\dot{\bar z}_i+\Gamma_i,\\
	\dot{\bar z}_i=-{\bf h}_i-\sum\nolimits_{j=1}^N a_{ij}(\bar \lm_i-\bar \lm_j),\\
	\dot{\bar \lm}_i=\sum\nolimits_{j=1}^N a_{ij}(\bar z_i-\bar z_j),\\
	\dot{\bar {\chi}}_i^a = \Pi \bar{\chi}_i^a +\Psi \bar{\chi}_i^b,\\
	\dot{\bar e}_i=(A_l\otimes I_m)\bar e_i\\
	\end{cases}
	$$
	where $\Gamma_i$ is Lipschitz in $(\bar \chi_i^a,\, \bar e_i)$ determined by \eqref{ctrl:decoupling} and \eqref{ctr:opt-offline-output-linear}.
	
	This whole system is again in a cascaded form with the last two subsystems as the driving one. By Theorem \ref{thm:opt-offline-linear} and Lemma \ref{lem:cascaded}, this implies the asymptotic stability when $f_i(\cdot)$ is strictly convex and global exponential stability under Assumption \ref{ass:convexity-strong}. Hence, the proof is complete.
\pe 

\begin{rem}
	Based on the above analysis,  it is worthwhile to mention that we can solve the problem by making the optimization design and tracking control design of high-order dynamics almost independent, because the generator for the optimization, independent of the high-order dynamics, is ``separately" designed by each agent, and then embedded into the designed reference-tracking control of the high-order dynamics. As a result, the design complexities brought by high-order dynamics are decoupled from those by the optimization task. By solving the two simpler subproblems, the optimal output consensus problem can be solved via constructive controllers.  Also, this embedded control framework enjoys a large flexibility in choosing optimal signal generators and tracking controllers, and therefore, may be useful in various optimization algorithms.
\end{rem}

\subsection{Case II}
In some cases, the function $\nabla f_i(\cdot)$ is hard to obtain and only the real-time gradient $\nabla f_i(y_i)$ is available, and then the optimal signal generator cannot be independently implemented. Suppose we substitute $\nabla f_i(z_i)$ by $\nabla f_i(y_i)$ in \eqref{sys:n=1}, there will be a mismatching error $\nabla f_i(z_i)-\nabla f_i(y_i)$ in the control.  As stated in Lemma \ref{lem:generator}, the proposed optimal signal generator under Assumption \ref{ass:convexity-strong} is actually robust with respect to additive perturbations.  Thus, these mismatching errors can be handled and compensated by its robustness along with small-gain techniques and a high-gain design.

We first show the result on high-order integrators for simplicity.

\begin{thm}\label{thm:opt-online}
	Suppose Assumptions \ref{ass:vrd},~\ref{ass:graph} and \ref{ass:convexity-strong} hold. Then there
	exists a constant $\epsilon^*>0$ such that the optimal output consensus problem \eqref{opt:main} with agents of the form
	$x^{(n)}_i=v_i,\; y_i=x_i$ is solved by the following control algorithm for any
	$\epsilon\in (0, \,\epsilon^*)$:
	\begin{align}\label{ctr:opt-online}
	\begin{cases}
	v_i=-\frac{1}{\e^n}[c_0(x_i-z_i)+\sum_{k=1}^{n-1}\e^{k} c_{k} x_i^{(k)}],\\
	\dot{z}_i=-\nabla f_i(y_i)-\sum\nolimits_{j=1}^N a_{ij}(\lm_i-\lm_j),\\
	\dot{\lm}_i=\sum\nolimits_{j=1}^N a_{ij}(z_i-z_j)
	\end{cases}
	\end{align}
	where $c_0,\, \dots,\, c_{n-1}$ are constants such that
	$\sum_{k=0}^{n-1}c_{k}s^{k}+s^n$ is Hurwitz.   Moreover, $y_i(t)$
	converges to the optimal solution $y^*$ exponentially as $t\to \infty$ for $i=1,\,\ldots,\,N$.
\end{thm}

\pb  Letting $\bar z_i=z_i-y^*$ and $\bar \lm_i=\lm_i-\lm_i^*$ gives
	$$
	\begin{cases}
	\dot{\hat x}_i= \frac{1}{\e} A_m\hat x_i+b_m[-{\bf h}_i-{\Delta}_i-\sum\nolimits_{j=1}^N a_{ij}(\bar \lm_i-\bar \lm_j)],\\
	\dot{\bar z}_i=-{\bf h}_i-{\Delta}_i-\sum\nolimits_{j=1}^N a_{ij}(\bar \lm_i-\bar \lm_j),\\
	\dot{\bar \lm}_i=\sum\nolimits_{j=1}^N a_{ij}(\bar z_i-\bar z_j)
	\end{cases}
	$$
	where ${\bf h}_i=\nabla f_i(z_i)-\nabla f_i(y^*)$ and ${\Delta}_i=\nabla f_i(y_i)-\nabla f_i(z_i)$,  or in a compact form
	\begin{align}
	\begin{cases}
	\dot{\hat x}= \frac{1}{\e} (I_N\otimes A_m) \hat x+(I_N\times b_m)[-{\bf h}-{ \Delta}-(L\otimes I_m)\bar \lm],\\
	\dot{\bar z}=-{\bf h}-{ \Delta}-(L\otimes I_m)\bar \lm,\\
	\dot{\bar \lm}=(L\otimes I_m)\bar z
	\end{cases}
	\end{align}
	where $\hat x=\mbox{col}(\hat x_1,\,\dots,\,\hat x_N)\in \R^{nN},\,{\bf h}$ and ${\Delta}$ are determined by ${\bf h}_i $ and ${\Delta_i}$.
	
	Perform another coordinate transformation: $\hat \lm_1=(r^\top
	\otimes I_m)\bar \lm,\, \hat \lm_2=(R^\top  \otimes I_m) \bar \lm$.
	Since $(r^\top \otimes I_m)\dot{\bar \lm}={\bf 0}$, it follows that
	$\hat \lm_1\equiv {\bf 0}$ and
	\begin{align}\label{sys:high-gain-thm-reduced}
	\begin{cases}
	\dot{\hat x}= \frac{1}{\e} (I_N\otimes A_m) \hat x+(I_N\times b_m)\bar \Gamma,\\
	\dot{\bar z}=-{\bf h}-{\Delta}-(LR\otimes I_m)\hat \lm_2,\\
	\dot{\hat \lm}_2=(R^\top L\otimes I_m)\bar z
	\end{cases}
	\end{align}
	where $\bar \Gamma=-{\bf h}-{\Delta}-(LR\otimes I_m)\hat \lm_2$.
	
	Let $\hat z\triangleq\mbox{col}(\bar z,\,\hat \lm_2)$.  It can be easily verified that ${\Delta}$ is $\vartheta$-Lipschitz in $x-z$ and hence $\hat x$. Note that $\bar \Gamma$ is also $\bar \vartheta_1$-Lipschitz in $\hat \lm_2$ for a positive constant $\bar \vartheta_1$. Thus, there exists positive constants $\bar \vartheta_2$ and $\bar \vartheta_3$ such that $||\bar \Gamma||\leq \bar \vartheta_2 ||\hat x||+ \bar \vartheta_3||\hat z||$.
	
	Next, we invoke a small-gain technique on \eqref{sys:high-gain-thm-reduced} to prove this theorem via tuning
	$\e$. By Lemma \ref{lem:key-stability} or the proof of Lemma
	\ref{lem:generator}, the following (nominal) system
	\begin{align}
	\begin{cases}
	\dot{\bar z}=-{\bf h}-(LR\otimes I_m)\hat \lm_2,\\
	\dot{\bar \lm}_2=(R^\top L\otimes I_m)\bar z
	\end{cases}
	\end{align}
	is globally exponentially stable under the given assumptions.
	
	Recalling the Lipschitzness of ${\bf h}$ in $\bar z$, we then apply
	the converse Lyapunov theorem (Theorem 4.15 in
	\cite{khalil2002nonlinear}) to this system, that is, there is a
	continuously differentiable Lyapunov function $V_1(\cdot)$ such that
	\begin{eqnarray*}
		&\hat c_1||\hat z||^2\leq V_1(\hat z)\leq \hat c_2 ||\hat z||^2,&\\
		&\frac{\partial V_1}{\partial \bar z}[-{\bf h}-(LR\otimes I_m)\hat \lm_2]+\frac{\partial V_1}{\partial \bar \lm_2}[(R^\top L\otimes I_m)z]\leq -\hat c_3 ||\hat z||^2,&\\
		&||\frac{\partial V_1}{\partial \hat z}||\leq \hat c_4||\hat z||&
	\end{eqnarray*}
	for some positive constants $\hat c_1,\,\hat c_2,\,\hat c_3$ and $\hat c_4$.
	
	Since $A_m$ is Hurwitz, there is a positive definite matrix $P_m\in
	\R^{nm\times nm}$ satisfying $A_m^\top P_m+P_mA_m=-I_{nm}$. We take
	a quadratic Lyapunov function $V=\hat x^\top (I_N\otimes P_m) \hat
	x+V_1(\hat z)$, whose derivative along the trajectory of system
	\eqref{sys:high-gain-thm-reduced} is
	\begin{align*}
	\dot{V}=&2\hat x^\top (I_N\otimes P_m)[\frac{1}{\e} (I_N\otimes A_m) \hat x+(I_N\times b_m)\bar \Gamma-\frac{\partial V_1}{\partial \hat z}\dot{\hat z}\\
	\leq &\frac{1}{\e}\hat x^\top [I_N\otimes (P_mA_m+A_m^\top P_m)]\hat x+2\hat x^\top (I_N\otimes P_mb_m)\bar \Gamma\\
	&+\frac{\partial V_1}{\partial \bar z}[-{\bf h}-(LR\otimes I_m)\hat \lm_2]+\frac{\partial V_1}{\partial \bar \lm_2}(R^\top L\otimes I_m)z-\frac{\partial V_1}{\partial \bar z}\Delta\\
	\leq &-\frac{1}{\e}||\hat x||^2+\hat \vartheta_2||\hat x||^2+\hat \vartheta_3||\hat x||||\hat z||-\hat c_3||\hat z||^2+\hat c_4 \vartheta||\hat x||||\bar z||\\
	\leq &-(\frac{1}{\e}-\hat \vartheta_2)||\hat x||^2+\hat \vartheta_4||\hat x||||\hat z||-\hat c_3||\hat z||^2
	\end{align*}
	where $\hat \vartheta_2=2\bar \vartheta_2||(I_N\otimes P_mb_m)||$, $\hat \vartheta_3=2\bar
	\vartheta_3||(I_N\otimes P_mb_m)||$, and $\hat \vartheta_4=\hat \vartheta_3+\hat c_4 \vartheta$.
	
	By the Young's inequality, we have
	\begin{align*}
	\dot{V}&\leq -(\frac{1}{\e}-\hat \vartheta_2)||\hat x||^2+\frac{\hat \vartheta_4^2}{2\hat c_3}||\hat x||^2+\frac{\hat c_3}{2}||\hat z||^2-\hat c_3||\hat z||^2\\
	&=-(\frac{1}{\e}-\hat \vartheta_2-\frac{\hat \vartheta_4^2}{2\hat c_3})||\hat x||^2-\frac{\hat c_3}{2}||\hat z||^2.
	\end{align*}
	We take
	\begin{align}\label{para:eps}
	\e^*=\frac{2\hat c_3}{2\hat c_3\hat \vartheta_2+ \hat \vartheta_4^2}
	\end{align}
	and obtain that for any $\e\in (0,\,\e^*)$
	\begin{align*}
	\dot{V}&\leq -(\frac{1}{\e}-\frac{1}{\e^*})||\hat x||^2-\frac{\hat c_3}{2}||\hat z||^2
	\end{align*}
	which implies the exponential stability of
	\eqref{sys:high-gain-thm-reduced} and therefore, the exponential
	convergence of $x_i$ with respect to $y^*$ as $t\to \infty$. Thus,
	the proof is complete. \pe 

Based on Theorem \ref{thm:opt-online} and the normal form
\eqref{sys:normal-form}, we go back to investigate linear
minimum-phase agents of the form \eqref{sys:agent} and the following
controller is readily obtained, that is,
\begin{equation}\label{ctr:linear-online}
\begin{cases}
\tilde{u}_i=\tilde R^{-1} (-\Upsilon \chi_i^a -S \chi_i^b +u_i),\\
v_i=-\frac{1}{\e^n}[c_0(y_i-z_i)+\sum_{k=1}^{n-1}\e^{k} c_{k} y_i^{(k)}],\\
\dot{z}_i=-\nabla f_i(y_i)-\sum\nolimits_{j=1}^N a_{ij}(\lm_i-\lm_j),\\
\dot{\lm}_i=\sum\nolimits_{j=1}^N a_{ij}(z_i-z_j)
\end{cases}
\end{equation}
where $c_1, \dots, c_{n-1}$ are selected as that in Theorem
\ref{thm:opt-online}, and the relationship between $u_i$ and $v_i$
is described by \eqref{ctrl:precompensator}.

Here is our main result of \eqref{ctr:linear-online} on general linear systems based on the embedded control scheme.

\begin{thm}\label{thm:opt-online-linear}
	Suppose Assumptions \ref{ass:vrd},~\ref{ass:graph} and \ref{ass:convexity-strong} hold. Then there
	is a constant $\epsilon^*>0$ such that the optimal output consensus problem \eqref{opt:main} with agents of the form
	\eqref{sys:agent} is solved by the control \eqref{ctr:linear-online}
	for any $\epsilon\in (0, \,\epsilon^*)$.  Moreover, $y_i(t)$
	converges to $y^*$ exponentially as $t\to \infty$ for $i=1,\,\ldots,\,N$.
\end{thm}

\pb 
	Following similar arguments in the proof of Theorem \ref{thm:opt-offline-linear}.  Clearly, the subsystem \eqref{ssys}-\eqref{output} is equivalent to
	the $n$th-order integrator.  Let ${\tilde
		\chi}_i^a={\chi}_i^a-{{\chi}_i^a}^*, \,{\tilde
		\chi}_i^b={\chi}_i^b-{{\chi}_i^b}^*$ (where $({{\chi}_i^a}^*,
	\,{{\chi}_i^b}^*)$ is the equilibrium point of
	\eqref{zds}-\eqref{ssys} determined by $y^*$) and take a similar
	coordinate transformation as that given in Theorem
	\ref{thm:opt-online}. Then the whole error system can be
	expressed in a cascaded form:
	$$
	\begin{cases}
	\dot{\tilde \chi}_i^a = \Pi \tilde\chi_i^a +\Psi \tilde\chi_i^b,\\
	\dot{\hat \chi}_i^c= \frac{1}{\e} A_m {\hat \chi}_i^c+b_m[-{\bf h}_i-{\Delta}_i-\sum\nolimits_{j=1}^N a_{ij}(\bar \lm_i-\bar \lm_j)],\\
	\dot{\bar z}_i=-{\bf h}_i-{\Delta}_i-\sum\nolimits_{j=1}^N a_{ij}(\bar \lm_i-\bar \lm_j),\\
	\dot{\bar \lm}_i=\sum\nolimits_{j=1}^N a_{ij}(\bar z_i-\bar z_j).
	\end{cases}
	$$
	
	Note that the matrix $\Pi$ is Hurwitz and $\Psi \tilde\chi_i^b$ is Lipschitz with respect to $({\hat \chi}_i^c,\, \bar z_i)$. Set $\e^*$ as defined in \eqref{para:eps} and by the proof of Theorem \ref{thm:opt-online}, the subsystem $(\hat \chi_i^c,\,\bar z_i,\,\bar \lm_i)$ is then globally exponentially stable with respect to its equilibrium point under control input $v_i$ for any $0< \e < \e^*$. Global exponential stability of the above cascaded system can be obtained by Lemma \ref{lem:cascaded}, which implies that the state trajectory of this agent is bounded and its output $y_i$ converges to $y^*$ exponentially as $t\to \infty$. The proof is thus complete.
\pe 

It can be found that for the given multi-agent system, the tracking control of each agent can be achieved by tuning only one input parameter $\e$. For different optimization problems, we can simply replace the local cost function $f_i(\cdot)$ in the optimal signal generator and then adjust the gain parameter $\e$, which transforms a controller design problem into a parameter-tuning one. Therefore, the embedded control scheme is much simpler than redesigning a complete new algorithm for these high-order agents, which may be favorable in large-scale networks.

\begin{rem}
	Similar problems have been investigated for the first-order system \cite{shi2013reaching, kia2015distributed} with vector relative degree $(1,\,...,\,1)$ and the second-order system \cite{zhang2014distributed, Deng} with vector relative degree $(2,\,...,\,2)$.  It is remarkable that the optimal signal generator used here actually has decoupled the complexity of optimization task from that of the high-order dynamics tracking problem, which may facilitate the design for even more
	complex agent dynamics.
\end{rem}

Next, let us consider its output feedback version.  As in Case I, by attaching the high-gain observer \eqref{eq.mpobserver}, we have the following distributed control for system \eqref{sys:agent}.
\begin{equation}\label{ctr:opt-online-output-linear}
\begin{cases}
\tilde{u}_i=\tilde R^{-1} (-\Upsilon \hat \chi_i^a -S \hat\chi_i^b +u_i),\\
v_i=-\frac{1}{\e^n}[c_0(y_i-z_i)+\sum_{k=1}^{n-1}\e^{k} c_{k} \hat y_i^{(k)}],\\
\dot{z}_i=-\nabla f_i(y_i)-\sum\nolimits_{j=1}^N a_{ij}(\lm_i-\lm_j),\\
\dot{\lm}_i=\sum\nolimits_{j=1}^N a_{ij}(z_i-z_j)\\
\end{cases}
\end{equation}
where $c_0, \, \ldots, \, c_{n-1}$ are defined as before and $\e>0$ is to be determined later.

We then have the following theorem for this design.

\begin{thm}\label{thm:online-ouput-linear}
	Suppose Assumptions \ref{ass:vrd},~\ref{ass:graph} and \ref{ass:convexity-strong} hold. Then there
	is a constant $\epsilon^*>0$ such that the optimal output consensus problem \eqref{opt:main} with agents of the form
	\eqref{sys:agent} is solved by \eqref{ctr:opt-online-output-linear}
	for any $\epsilon\in (0, \,\epsilon^*)$. Moreover, $y_i(t)$ converge
	to the optimal solution $y^*$ exponentially as $t\to \infty$ for $i=1,\,\ldots,\,N$.
\end{thm}

\pb   The proof is similar to Theorem \ref{thm:offline-ouput-linear}. In fact, following the same arguments, we only have to prove the stability of the whole error system: 
	$$
	\begin{cases}
	\dot{\tilde \chi}_i^a = \Pi \tilde\chi_i^a +\Psi \tilde\chi_i^b,\\
	\dot{\hat \chi}_i^c= \frac{1}{\e} A_m {\hat \chi}_i^c+b_m[-{\bf h}_i-{\Delta}_i-\sum\nolimits_{j=1}^N a_{ij}(\bar \lm_i-\bar \lm_j)]+\Gamma_i,\\
	\dot{\bar z}_i=-{\bf h}_i-{\Delta}_i-\sum\nolimits_{j=1}^N a_{ij}(\bar \lm_i-\bar \lm_j),\\
	\dot{\bar \lm}_i=\sum\nolimits_{j=1}^N a_{ij}(\bar z_i-\bar z_j),\\
	\dot{\bar {\chi}}_i^a = \Pi \bar{\chi}_i^a +\Psi \bar{\chi}_i^b,\\
	\dot{\bar e}_i=(A_l\otimes I_m)\bar e_i\\
	\end{cases}
	$$
	where $\Gamma_i$ is Lipschitz in $(\bar \chi_i^a,\, \bar e_i)$ determined by \eqref{ctrl:decoupling} and \eqref{ctr:opt-online-output-linear}.
	
	This system is again in a cascaded form. Choosing $\e^*$ as defined in \eqref{para:eps}, one can obtain its global exponential stability for any $\e\in (0,\,\e^*)$ under Assumption \ref{ass:convexity-strong} by Lemma \ref{lem:cascaded} and Theorem \ref{thm:opt-online-linear}. Hence, the proof is complete.
\pe 

\begin{rem}\label{rem:robustness}
	In practice,  bounded uncertainties may come from the communication noise and gradient calculation.  Suppose Assumptions \ref{ass:vrd}, \ref{ass:graph} and \ref{ass:convexity-strong} hold, it can be easily verified that our algorithm is robust in the sense of bounded-input bounded-output stability with respect to those uncertainties as its input and approximation error of $y_i-y_i^*$ as its output. In fact, we can make the approximation error arbitrarily small by adjusting the gain parameter $\e$.
\end{rem}

\begin{rem}
	As we remarked before, the optimal output consensus problem is a combination of the two hot topics: output consensus of high-order agents\cite{ren2007distributed, ma2010necessary, xi2012output, seo2009consensus} and distributed optimization for single integrators \cite{kia2015distributed, lin2016distributed, yi2014quantized, shi2013reaching}. Due to the couplings between the distributed optimization requirement and high-order dynamic processes, it is much more challenging than traditional distributed optimization or output consensus problem. This is why we come up with the optimal signal generator to decouple the optimization design from the dynamics in order to simplify the whole design complexities. Compared with the new result given in \cite{zhang2015distributed},  this embedded scheme leads to a generally constructive way to solve the problem.
\end{rem}

\section{Simulations}\label{sec:simu}

In this section, we present two examples to illustrate our problem and the effectiveness of our designs.

{\it Example 1}. Consider an optimal rendezvous problem \cite{ren2007distributed} of wheeled robots with the following dynamics.
\begin{align*}
\begin{cases}
\dot{r}_i^x=v_i\cos(\theta_i),\\
\dot{r}_i^y=v_i\sin(\theta_i),\\
\dot{\theta}_i=\omega_i,\\
\dot{v}_i=\frac{1}{m_i}F_i,\\
\dot{\omega}_i=\frac{1}{J_i}\tau_i
\end{cases}
\end{align*}
where $(r_i^x,\, r_i^y, \,\theta_i)$ are the inertia center's position and orientation of the $i$th robot, $(v_i,\,\omega_i)$ the linear and angular speed, $(F_i,\,\tau_i)$ the applied force and torque, $(m_i,\, J_i)$ the mass and moment of inertia for $i=1,\,\dots,\,5$. Let $d_i$ represent the distance between the hand position and inertia center of the $i$th robot. Applying feedback linearization as in \cite{ren2007distributed} about the hand position $\tilde x_i \triangleq (r_i^x+d_i\cos(\theta_i), \,r_i^y+d_i\sin(\theta_i))$ yields a simple linear dynamics as follows:
\begin{align*}
\dot{\tilde x}_i=\tilde v_i,\, \dot{\tilde  v}_i=\tilde u_i,\, y_i=\tilde x_i.
\end{align*}
It apparently satisfies Assumption 1 with relative degree $(2,\,2)$.

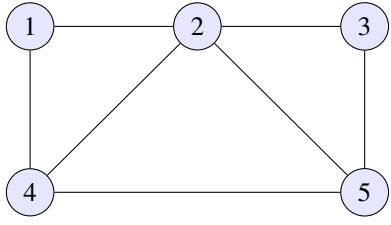
\begin{figure}
	\centering
	\begin{tikzpicture}[node distance=2.2 cm, >=stealth',
	every state/.style ={circle, fill=blue!10, minimum size=0.10cm}]
	\node[align=center,state](node1) {1};
	\node[align=center,state](node2)[right of=node1]{2};
	\node[align=center,state](node3)[right of=node2]{3};
	\node[align=center,state](node4)[below of=node1]{4};
	\node[align=center,state](node5)[below of=node3]{5};
	\path[-]    (node1) edge (node2)
	(node2) edge  (node3)
	(node3) edge  (node5)
	(node1) edge  (node4)
	(node2) edge  (node4)
	(node2) edge  (node5)
	(node4) edge  (node5)  ;
	\end{tikzpicture}\\
	\caption{Interaction topology of the multi-agent system.}\label{fig:graph}
\end{figure}

To drive all hands of robots to rendezvous at a common point that minimizes the aggregate distance from their starting points to this final location, the cost functions satisfying Assumption \ref{ass:convexity-strong} are as $f_i(y)=\frac{1}{2} ||y-y_i(0)||^2$ and $f(y)=\frac{1}{2}\sum_{i=1}^5 ||y-y_i(0)||^2$ ($i=1,\,\dots,\, 5$). We can easily check that the optimal solution of the global cost function is $\mbox{Aver(y(0))}\triangleq \frac{1}{N}\sum_{i=1}^Ny_i(0)$. The communication graph among these robots satisfying Assumption \ref{ass:graph} is depicted as Fig.~\ref{fig:graph} with all the edge weights as $1$. Recalling Theorem \ref{thm:opt-offline}, the rendezvous problem can be solved by controller \eqref{ctr:opt-offline}.

Because $\nabla f_i(y)=y-y_i(0)$, the generator \eqref{sys:n=1} reduces to:
\begin{eqnarray*}
	\begin{cases}
		\dot{z}_i=-(z_i-y_i(0))-\sum\nolimits_{j=1}^N a_{ij}(\lm_i-\lm_j),\\
		\dot{\lm}_i=\sum\nolimits_{j=1}^N a_{ij}(z_i-z_j).
	\end{cases}
\end{eqnarray*}
Thus, the state feedback control \eqref{ctr:opt-offline} is as follows:
\begin{eqnarray*}
	\begin{cases}
		\tilde{u}_i=-\frac{1}{\e^2}[c_0(y_i-z_i)+\e c_{1} \tilde v_i],\\
		\dot{z}_i=-(z_i-y_i(0))-\sum\nolimits_{j=1}^N a_{ij}(\lm_i-\lm_j),\\
		\dot{\lm}_i=\sum\nolimits_{j=1}^N a_{ij}(z_i-z_j)
	\end{cases}
\end{eqnarray*}

\begin{figure}
	\centering
	\includegraphics[width=0.75\textwidth]{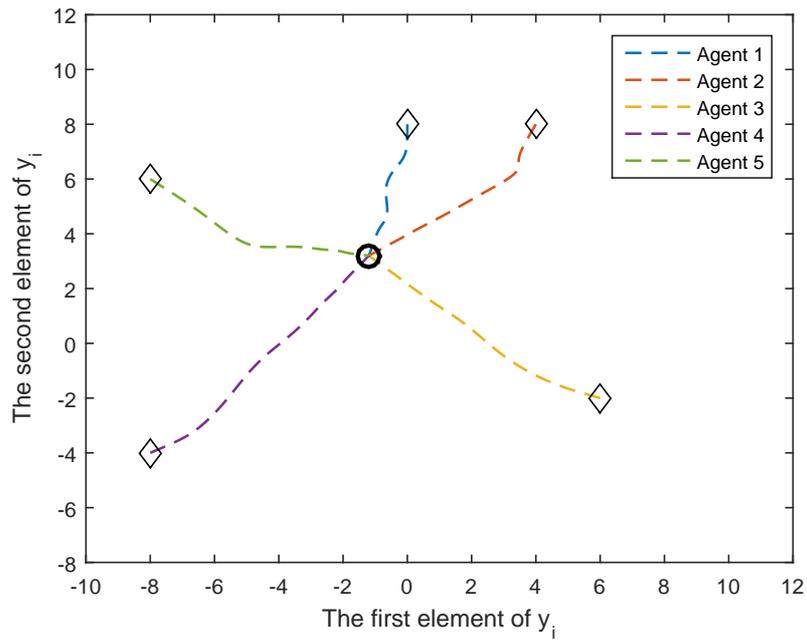}
	\caption{Rendezvous of five mobile robots under state feedback control~\eqref{ctr:opt-offline}. }\label{fig:simu-aver}
\end{figure}

Take $c_0=4, \, c_1=8,\, \e=1$ and all initials (randomly) in $[-10,\,10]^{8}$. The local minimizers, i.e. their initial points, are marked by diamonds, the global optimal solution $y^*$ by a circle. The simulation result is given in Fig.~\ref{fig:simu-aver} and all robots achieve the optimal rendezvous at $y^*$.

{\it Example 2}. To verify the effectiveness of our embedded design, we then consider a high-order multi-agent system modified from \cite{xi2012output} with more complex objective functions.
\begin{eqnarray}\label{sys:simsys}
\dot{\tilde{x}}_i=\tilde{A} {\tilde x}_i + \tilde B \tilde u_i, \quad y_i =\tilde C \tilde{x}_i,
\quad i=1,\,\ldots, \,5
\end{eqnarray}
where
$$
\tilde A =\begin{bmatrix}   0.1  &0.1  &-2.5&0.5\\
0  &-0.4&0&-1.5\\
2  &0  &0  &-1\\
0  &0.2&1  &0.8
\end{bmatrix},~~ \tilde B = \begin{bmatrix}
1&0\\
0&1\\
0&0\\
1&1
\end{bmatrix},$$
$$\tilde C = \begin{bmatrix}
0&0&1&0\\
0&0&0&1
\end{bmatrix}.
$$

The local cost functions are as follows.
\begin{align*}
{f_1}(y) =& (y_a - 8)^2+(y_b-1)^2,\\
{f_2}(y) =& \frac{y_a^2}{{20\sqrt {{{({y_a})}^2} + 1} }} + \frac{y_b^2}{{20\sqrt {{{({y_b})}^2} + 1} }} + ||y ||^2,\\
{f_3}(y) =& \frac{{{{({y_a})}^2}}}{{80\ln \left( {{{({y_a})}^2} + 2} \right)}} + \frac{{{{({y_b})}^2}}}{{80\ln \left( {{{({y_b})}^2} + 2} \right)}}
+ {\left\| {y - 5 \times {\bf 1}_2} \right\|^2},\\
{f_4}(y) =& \ln \left( {{e^{ - 0.05{y_a}}} + {e^{0.05{y_a}}}} \right) + \ln \left( {{e^{ - 0.05{y_b}}} + {e^{0.05{y_b}}}} \right)
+ {\left\| y \right\|^2},\\
{f_5}(y) =& \frac{y_a^2}{{25\sqrt {{{({y_a})}^2} + 1} }} + \frac{y_b^2}{{25\sqrt {{{({y_b})}^2} + 1} }} + {\left\| y \right\|^2} + {\bf 1}_2^\top y
\end{align*}
where $y=[y_{a} ~ y_{b}]^T$. The interaction topology of this multi-agent system is still taken as Fig.~\ref{fig:graph}. It can be verified that the vector relative degree of agents is $(2,\,1)$ and Assumptions \ref{ass:vrd}-- \ref{ass:convexity-strong} hold. Additionally, $y_a^*=2.5, \,y_b^*=1.1$ by minimizing the global cost function $f(y)=\sum_{i=1}^5 f_i^*(y)$.

Assume that agents only have real-time gradients. Recalling Theorem \ref{thm:online-ouput-linear}, this output optimal consensus problem can be solved by output feedback controller \eqref{ctr:opt-online-output-linear}.

In fact, by taking $\phi_{\iota}=[\,1,\,1,\,0,\,-1\,]$, the system \eqref{sys:simsys} can be put into \eqref{sys:normal-form-vrd} as follows:
\begin{eqnarray*}
	\left\{\begin{split}
		\dot{\chi}_i^a &= -0.5 \chi_i^a +[-3.5,\,0.3,\,-2] \chi_i^b, \\
		\dot{\xi}^{ib}_{1 1} &= \xi^{ib}_{1 2},\\
		\dot{\xi}^{ib}_{1 2} &= [-6,\,0.1,\,0.3] \chi_i^b + [1,\,-1] \tilde{u}_i, \\
		\dot{\xi}^{ib}_{2 1} &= 0.2 \chi_i^a +[1,\,-0.1,\,0.9] \chi_i^b + [1,\,1] \tilde{u}_i,\\
		y_i &=\mbox{col}(\xi^{ib}_{11}, \xi^{ib}_{21}).
	\end{split}\right.
\end{eqnarray*}
where $\tilde R=[1,\,-1;\, 1,\,1]$.

Thus, the whole controller is given as follows:
\begin{eqnarray*}
	\begin{cases}
		\tilde{u}_i=\begin{bmatrix}
			-0.1 \hat{\chi}_i^a + 2.5 \hat{\xi}_{11}^{ib}-0.6 \hat{\xi}_{21}^{ib} +0.5 u_{i1}+0.5 u_{i2} \\
			-0.1 \hat{\chi}_i^a - 3.5 \hat{\xi}_{11}^{ib}-0.3 \hat{\xi}_{21}^{ib} +0.1\hat{\xi}_{12}^{ib}-0.5 u_{i1}+0.5 u_{i2} \\
		\end{bmatrix},\\
		u_{i1} =v_{i1}, ~ \dot{u}_{i2} =v_{i2}, \\
		v_i = \begin{bmatrix}
			v_{i1} \\
			v_{i2} \\
		\end{bmatrix}=-\frac{1}{\e^2}[c_0(y_i-z_i)+\e c_{1} y_i^{(1)}],\\
		\dot{z}_i=-\nabla f_i(y_i)-\sum\nolimits_{j=1}^N a_{ij}(\lm_i-\lm_j),\\
		\dot{\lm}_i=\sum\nolimits_{j=1}^N a_{ij}(z_i-z_j)
	\end{cases}
\end{eqnarray*}
where
$\hat{y}_i=\mbox{col}(\hat{\xi}_{11}^{ib},\,\hat{\xi}_{21}^{ib}),\,
\hat{y}_i^{(1)}= \mbox{col}(\hat{\xi}_{12}^{ib},\,\hat{u}_{i2})$,
and $\hat{\chi}_i^a$ are generated by the following observer
\begin{eqnarray*}
	\left\{\begin{split}
		\dot{\hat{y}}_i&=\hat{y}_i^{(1)} - l_1 (\hat{y}_i-y_i), \\
		\dot{\hat{y}}_i^{(1)}&=v_i-  l_2(\hat{y}_i-y_i),  \\
		\dot{\hat{\chi}}_i^a &=  -3 \hat{\chi}_i^a + [-3.5,\,0.3,\,-2] \hat {\chi}_i^b.
	\end{split}\right.
\end{eqnarray*}

\begin{figure}
	\centering
	\includegraphics[width=0.75\textwidth]{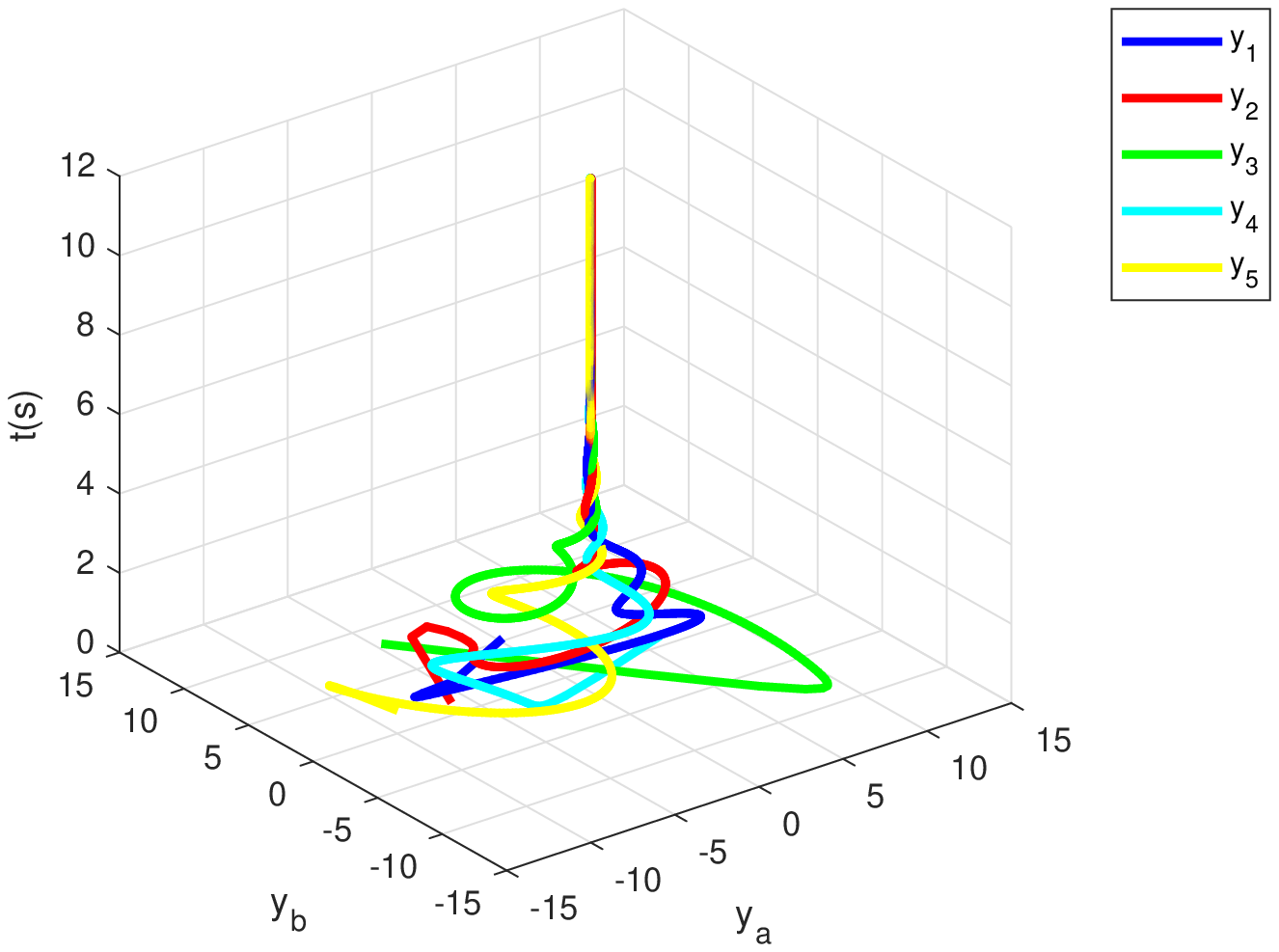}
	\caption{Phase portraits of $y_i$ under output feedback control \eqref{ctr:opt-online-output-linear}.}\label{fig:simu-1}
\end{figure}

\begin{figure}
	\centering
	\includegraphics[width=0.75\textwidth]{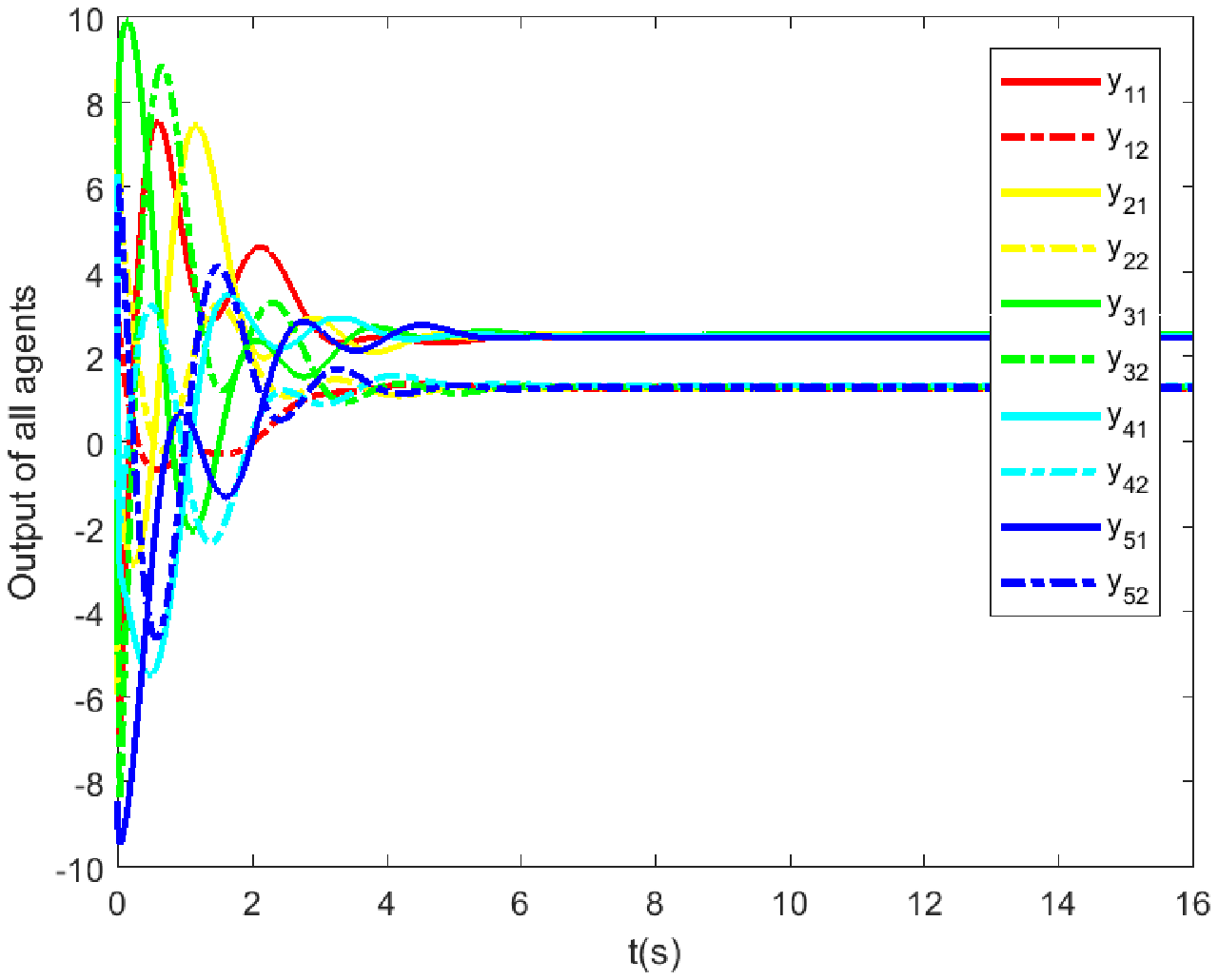}
	\caption{Profiles of all components in $y_i$ under output feedback control \eqref{ctr:opt-online-output-linear}.}\label{fig:simu-2}
\end{figure}

Take $c_0=4,\, c_1=8,\, l_1=4,\, l_2=8,\, \e=1$, and all initials are generated (randomly) in $[-10,\,10]^{14}$. This output optimal consensus problem can be solved by the output feedback controller \eqref{ctr:opt-online-output-linear}.  The simulation results are given in Figs.~\ref{fig:simu-1} and \ref{fig:simu-2} and all the outputs converge to the global optimal point.

\section{Conclusions}\label{sec:con}
This paper has investigated the optimal output consensus problem for high-order minimum-phase multi-agent systems.  An embedded control scheme has been proposed and applied to solve this problem based on the introduction of an optimal signal generator. The proposed algorithms have been proved to converge to the optimal solution asymptotically or exponentially with different conditions. In fact, many challenging optimal output consensus problems remain to be done, including the cases of practical and nonlinear agents, or various uncertainties from communication or
environment, or different optimization constraints.



\bibliographystyle{IEEEtran}
\bibliography{opt-high}

\begin{thebibliography}{10}
\providecommand{\url}[1]{#1}
\csname url@samestyle\endcsname
\providecommand{\newblock}{\relax}
\providecommand{\bibinfo}[2]{#2}
\providecommand{\BIBentrySTDinterwordspacing}{\spaceskip=0pt\relax}
\providecommand{\BIBentryALTinterwordstretchfactor}{4}
\providecommand{\BIBentryALTinterwordspacing}{\spaceskip=\fontdimen2\font plus
\BIBentryALTinterwordstretchfactor\fontdimen3\font minus
  \fontdimen4\font\relax}
\providecommand{\BIBforeignlanguage}[2]{{%
\expandafter\ifx\csname l@#1\endcsname\relax
\typeout{** WARNING: IEEEtran.bst: No hyphenation pattern has been}%
\typeout{** loaded for the language `#1'. Using the pattern for}%
\typeout{** the default language instead.}%
\else
\language=\csname l@#1\endcsname
\fi
#2}}
\providecommand{\BIBdecl}{\relax}
\BIBdecl

\bibitem{fax2004information}
J.~A. Fax and R.~M. Murray, ``Information flow and cooperative control of
  vehicle formations,'' \emph{IEEE Transactions on Automatic Control}, vol.~49,
  no.~9, pp. 1465--1476, 2004.

\bibitem{wang2012leader}
X.~Wang, W.~Ni, and X.~Wang, ``Leader-following formation of switching
  multirobot systems via internal model,'' \emph{IEEE Transactions on Systems,
  Man, and Cybernetics, Part B: Cybernetics}, vol.~42, no.~3, pp. 817--826,
  2012.

\bibitem{su2012cyber}
Y.~Su and J.~Huang, ``Cooperative output regulation with application to
  multi-agent consensus under switching network,'' \emph{IEEE Transactions on
  Systems Man and Cybernetics Part B: Cybernetics}, vol.~42, no.~3, pp.
  864--875, Jun. 2012.

\bibitem{tang2015distributed}
Y.~Tang, Y.~Hong, and X.~Wang, ``Distributed output regulation for a class of
  nonlinear multi-agent systems with unknown-input leaders,''
  \emph{Automatica}, vol.~62, pp. 154--160, 2015.

\bibitem{rezaee2015average}
H.~Rezaee and F.~Abdollahi, ``Average consensus over high-order multiagent
  systems,'' \emph{IEEE Transactions on Automatic Control}, vol.~60, no.~11,
  pp. 3047--3052, 2015.

\bibitem{seo2009consensus}
J.~H. Seo, H.~Shim, and J.~Back, ``Consensus of high-order linear systems using
  dynamic output feedback compensator: low gain approach,'' \emph{Automatica},
  vol.~45, no.~11, pp. 2659--2664, 2009.

\bibitem{xi2012output}
J.~Xi, Z.~Shi, and Y.~Zhong, ``Output consensus analysis and design for
  high-order linear swarm systems: partial stability method,''
  \emph{Automatica}, vol.~48, no.~9, pp. 2335--2343, 2012.

\bibitem{ren2007distributed}
W.~Ren and E.~Atkins, ``Distributed multi-vehicle coordinated control via local
  information exchange,'' \emph{International Journal of Robust and Nonlinear
  Control}, vol.~17, no. 10-11, pp. 1002--1033, 2007.

\bibitem{ma2010necessary}
C.~Ma and J.~Zhang, ``Necessary and sufficient conditions for consensusability
  of linear multi-agent systems,'' \emph{IEEE Transactions on Automatic
  Control}, vol.~55, no.~5, pp. 1263--1268, May 2010.

\bibitem{nedic2010constrained}
A.~Nedi{\'c}, A.~Ozdaglar, and P.~Parrilo, ``Constrained consensus and
  optimization in multi-agent networks,'' \emph{IEEE Transactions on Automatic
  Control}, vol.~55, no.~4, pp. 922--938, 2010.

\bibitem{boyd2011distributed}
S.~Boyd, N.~Parikh, E.~Chu, B.~Peleato, and J.~Eckstein, ``Distributed
  optimization and statistical learning via the alternating direction method of
  multipliers,'' \emph{Foundations and Trends in Machine Learning}, vol.~3,
  no.~1, pp. 1--122, 2011.

\bibitem{bose2015equivalent}
S.~Bose, S.~H. Low, T.~Teeraratkul, and B.~Hassibi, ``Equivalent relaxations of
  optimal power flow,'' \emph{IEEE Transactions on Automatic Control}, vol.~60,
  no.~3, pp. 729--742, 2015.

\bibitem{zhang2015twc}
Y.~Zhang, Y.~Lou, Y.~Hong, and L.~Xie, ``Distributed projection-based
  algorithms for source localization in wireless sensor networks,'' \emph{IEEE
  Transactions on Wireless Communications}, vol.~14, no.~6, pp. 3131--3142,
  Nov. 2015.

\bibitem{yuan2011distributed}
D.~Yuan, S.~Xu, and H.~Zhao, ``Distributed primal-dual subgradient method for
  multiagent optimization via consensus algorithms,'' \emph{IEEE Transactions
  on Systems, Man, and Cybernetics, Part B: Cybernetics}, vol.~41, no.~6, pp.
  1715--1724, 2011.

\bibitem{yi2014quantized}
P.~Yi and Y.~Hong, ``Quantized subgradient algorithm and data-rate analysis for
  distributed optimization,'' \emph{IEEE Transactions on Control of Network
  Systems}, vol.~1, no.~4, pp. 380--392, 2014.

\bibitem{lin2016distributed}
P.~Lin, W.~Ren, and Y.~Song, ``Distributed multi-agent optimization subject to
  nonidentical constraints and communication delays,'' \emph{Automatica},
  vol.~65, pp. 120--131, 2016.

\bibitem{shi2013reaching}
G.~Shi, K.~H. Johansson, and Y.~Hong, ``Reaching an optimal consensus:
  dynamical systems that compute intersections of convex sets,'' \emph{IEEE
  Transactions on Automatic Control}, vol.~58, no.~3, pp. 610--622, 2013.

\bibitem{kia2015distributed}
S.~S. Kia, J.~Cort{\'e}s, and S.~Mart{\'\i}nez, ``Distributed convex
  optimization via continuous-time coordination algorithms with discrete-time
  communication,'' \emph{Automatica}, vol.~55, pp. 254--264, 2015.

\bibitem{wang2015cyber}
X.~Wang, Y.~Hong, and H.~Ji, ``Distributed optimization for a class of
  nonlinear multiagent systems with disturbance rejection,'' \emph{IEEE
  Transactions on Cybernetics}, vol.~46, no.~7, pp. 1655--1666, July 2016.

\bibitem{liu2015projection}
Q.~Liu and J.~Wang, ``A projection neural network for constrained quadratic
  minimax optimization,'' \emph{IEEE Transactions on Neural Networks and
  Learning Systems}, vol.~26, no.~11, pp. 2891--2900, Nov. 2015.

\bibitem{yi2015distributed}
P.~Yi, Y.~Hong, and F.~Liu, ``Distributed gradient algorithm for constrained
  optimization with application to load sharing in power systems,''
  \emph{Systems \& Control Letters}, vol.~83, pp. 45--52, 2015.

\bibitem{liu2015second}
Q.~Liu and J.~Wang, ``A second-order multi-agent network for bound-constrained
  distributed optimization,'' \emph{IEEE Transactions on Automatic Control},
  vol.~60, no.~12, pp. 3310--3315, Dec. 2015.

\bibitem{qiu2016distributed}
Z.~Qiu, S.~Liu, and L.~Xie, ``Distributed constrained optimal consensus of
  multi-agent systems,'' \emph{Automatica}, vol.~68, pp. 209--215, 2016.

\bibitem{yang2016PI}
S.~Yang, Q.~Liu, and J.~Wang, ``A multi-agent system with a
  proportional-integral protocol for distributed constrained optimization,''
  \emph{IEEE Transactions on Automatic Control}, vol.~PP, no.~99, pp. 1--1,
  2016.

\bibitem{kim2014cooperative}
C.~Kim, D.~Song, Y.~Xu, J.~Yi, and X.~Wu, ``Cooperative search of multiple
  unknown transient radio sources using multiple paired mobile robots,''
  \emph{IEEE Transactions on Robotics}, vol.~30, no.~5, pp. 1161--1173, 2014.

\bibitem{anese2013distributed}
E.~Dall'Anese, H.~Zhu, and G.~B. Giannakis, ``Distributed optimal power flow
  for smart microgrids,'' \emph{IEEE Transactions on Smart Grid}, vol.~4,
  no.~3, pp. 1464--1475, Sep. 2013.

\bibitem{zhang2014distributed}
Y.~Zhang and Y.~Hong, ``Distributed optimization design for second-order
  multi-agent systems,'' in \emph{Control Conference (CCC), 2014 33rd
  Chinese}.\hskip 1em plus 0.5em minus 0.4em\relax IEEE, 2014, pp. 1755--1760.

\bibitem{Deng}
Z.~Deng and Y.~Hong, ``Multi-agent optimization design for autonomous
  lagrangian systems,'' \emph{Unmanned Systems}, vol.~4, no.~01, pp. 5--13,
  2016.

\bibitem{zhang2015distributed}
Y.~Zhang and Y.~Hong, ``Distributed optimization design for high-order
  multi-agent systems,'' in \emph{Control Conference (CCC), 2015 34th
  Chinese}.\hskip 1em plus 0.5em minus 0.4em\relax IEEE, 2015, pp. 7251--7256.

\bibitem{rockafellar1970convex}
R.~T. Rockafellar, \emph{Convex Analysis}.\hskip 1em plus 0.5em minus
  0.4em\relax Princeton, N.J.: Princeton University Press, 1970.

\bibitem{bertsekas2003convex}
D.~P. Bertsekas, A.~Nedic, and A.~E. Ozdaglar, \emph{Convex analysis and
  optimization}.\hskip 1em plus 0.5em minus 0.4em\relax Athena Scientific,
  2003.

\bibitem{godsil2001algebraic}
C.~Godsil and G.~F. Royle, \emph{Algebraic Graph Theory}.\hskip 1em plus 0.5em
  minus 0.4em\relax New York: Springer, 2001.

\bibitem{falb1967decoupling}
P.~L. Falb and W.~Wolovich, ``Decoupling in the design and synthesis of
  multivariable control systems,'' \emph{IEEE Transactions on Automatic
  Control}, vol.~12, no.~6, pp. 651--659, 1967.

\bibitem{isidori1995nonlinear}
A.~Isidori, \emph{Nonlinear Control Systems (3rd ed.)}, ser. Communications and
  control engineering series.\hskip 1em plus 0.5em minus 0.4em\relax Springer
  Science \& Business Media, 1995.

\bibitem{bhattacharya2011distributed}
S.~Bhattacharya and V.~Kumar, ``Distributed optimization with pairwise
  constraints and its application to multi-robot path planning,''
  \emph{Robotics: Science and Systems VI}, p. 177, 2011.

\bibitem{derenick2007convex}
J.~C. Derenick and J.~R. Spletzer, ``Convex optimization strategies for
  coordinating large-scale robot formations,'' \emph{IEEE Transactions on
  Robotics}, vol.~23, no.~6, pp. 1252--1259, 2007.

\bibitem{nesterov2013introductory}
Y.~Nesterov, \emph{Introductory Lectures on Convex Optimization: A Basic
  Course}.\hskip 1em plus 0.5em minus 0.4em\relax Springer Science \& Business
  Media, 2013, vol.~87.

\bibitem{hristu2005handbook}
D.~Hristu-Varsakelis and W.~S. Levine,
  \emph{\BIBforeignlanguage{English}{Handbook of {Networked} and {Embedded}
  {Control} {Systems}}}.\hskip 1em plus 0.5em minus 0.4em\relax Boston:
  Birkh{\"a}user, 2005.

\bibitem{chen1995linear}
C.-T. Chen, \emph{Linear system theory and design}.\hskip 1em plus 0.5em minus
  0.4em\relax Oxford University Press, Inc., 1995.

\bibitem{khalil2002nonlinear}
H.~K. Khalil, \emph{Nonlinear Systems (3rd ed.)}.\hskip 1em plus 0.5em minus
  0.4em\relax Upper Saddle River, N.J.: Prentice Hall, 2002.

\bibitem{ioannou1995robust}
P.~A. Ioannou and J.~Sun, \emph{Robust Adaptive Control}.\hskip 1em plus 0.5em
  minus 0.4em\relax Upper Saddle River, N.J.: Prentice-Hall, 1995.

\bibitem{alur2003hierarchical}
R.~Alur, T.~Dang, J.~Esposito, Y.~Hur, F.~Ivancic, V.~Kumar, P.~Mishra,
  G.~Pappas, and O.~Sokolsky, ``Hierarchical modeling and analysis of embedded
  systems,'' \emph{Proceedings of the IEEE}, vol.~91, no.~1, pp. 11--28, 2003.

\end{thebibliography}

\end{document}